\documentclass[a4paper,12pt]{amsart}
\usepackage{amssymb,amsmath,a4wide,graphicx,stmaryrd,fullpage,setspace,microtype,textcomp,tikz,tikz-cd,accents,mathtools,abraces}
\usepackage{hyperref}
\usepackage[shortlabels]{enumitem}
\definecolor{dark-red}{rgb}{0.5,0.15,0.15}
\hypersetup{
	colorlinks   = true, 
	urlcolor     = dark-red, 
	linkcolor    = black, 
	citecolor   = dark-red 
}
\usepackage[T1]{fontenc}
\usepackage{lmodern}
\usepackage[numbers,sort&compress]{natbib}

\title{Branching spaces of multipointed d-spaces}

\author[P. Gaucher]{Philippe Gaucher}

\address{Universit\'e Paris Cit\'e, CNRS, IRIF, F-75013, Paris, France}

\urladdr{http://www.irif.fr/{\~{}}gaucher}

\subjclass[2020]{55U35,55U99,68Q85}

\keywords{multipointed d-space, flow, branching, merging, globular subdivision, directed path, directed homotopy}

\swapnumbers

\makeatletter
\let\P\@undefined
\let\leq\@undefined
\let\geq\@undefined
\let\top\@undefined
\let\epsilon\@undefined
\let\vec\@undefined
\makeatother

\newcommand{\P}{\mathbb{P}}
\newcommand{\leq}{\leqslant}
\newcommand{\geq}{\geqslant}
\newcommand{\top}{{\mathbf{Top}}}
\newcommand{\epsilon}{\varepsilon}
\newcommand{\vec}{\overrightarrow}

\newcommand{\C}{\mathcal{C}}
\newcommand{\K}{\mathcal{K}}
\newcommand{\p}{\times}
\DeclareMathOperator{\hop}{ho\P}
\DeclareMathOperator{\pr}{pr}
\newcommand{\cont}{{\vec{\mathrm{Sp}}}}
\newcommand{\discont}{{\vec{\Omega}}}
\newcommand{\dtop}{{\brm{Flow}}}
\newcommand{\dcat}{{\mathrm{cat}}}
\newcommand{\glob}{{\rm{Glob}}}
\DeclareMathOperator{\cell}{{\brm{cell}}}
\newcommand{\sbd}{\mathrm{sbd}}
\newcommand{\iso}{\cong}
\newcommand{\vI}{\vec{I}}
\newcommand{\tr}[1]{{\langle{#1}\rangle}}
\DeclareMathOperator{\cl}{Cl}
\newcommand{\brm}[1]{\rm{\mathbf{#1}}}
\newcommand{\ttop}{{\brm{TOP}}}
\newcommand{\globM}{{\rm{Glob}}^{top}}
\DeclareMathOperator{\id}{Id}
\newcommand{\liminj}{\varinjlim}

\newcommand{\ptop}[1]{{\brm{{#1}dTop}}}
\newcommand{\cocartesian}{\arrow[lu, phantom, "\ulcorner"{font=\Large}, pos=0]}

\DeclareMathOperator{\im}{im}
\DeclareMathOperator{\dt}{dt}
\DeclareMathOperator{\natgl}{\underline{nat}^{gl}}

\newtheorem*{thmN}{Theorem}
\newtheorem{thm}{Theorem}[section]
\newtheorem{prop}[thm]{Proposition}
\newtheorem{lem}[thm]{Lemma}
\newtheorem{cor}[thm]{Corollary}
\newcommand{\bth}{\begin{thm}}
\renewcommand{\eth}{\end{thm}}
\newcommand{\bpf}{\begin{proof}}
\newcommand{\epf}{\end{proof}}
\theoremstyle{definition}
\newtheorem{defn}[thm]{Definition}
\newcommand{\bd}{\begin{defn}}
\newcommand{\ed}{\end{defn}}
\newtheorem{nota}[thm]{Notation}

\makeatletter
\newcommand*{\@opargbegintheorem}[3]{\trivlist
	\item[\hskip \labelsep{\bfseries #1\ #2}] \textbf{(#3)}\ \itshape}
\makeatother

\allowdisplaybreaks

\begin{document}

\begin{abstract} 
	Using the notion of a short directed path, we introduce the branching space of a multipointed $d$-space. We prove that for any q-cofibrant multipointed $d$-space, it is homeomorphic to the branching space of the q-cofibrant flow obtained by applying the categorization functor. As an application, we deduce a purely topological proof of the invariance of the branching space and of the branching homology of cellular multipointed $d$-spaces up to globular subdivision. By reversing the time direction, the same results are obtained for the merging space and the merging homology. 
\end{abstract}

\maketitle
\tableofcontents
\hypersetup{linkcolor = dark-red}

\section{Introduction}

\subsection*{Presentation}

Directed Algebraic Topology (DAT) provides a geometric framework for modeling concurrent systems in computer science, focusing on their behavior up to homotopy \cite{DAT_book}. At its core, DAT interprets an $n$-cube as a representation of the concurrent execution of $n$ actions, giving rise to precubical sets as natural models of concurrency \cite{zbMATH07226006}. The overarching objective is to construct algebraic invariants that remain stable under cubical subdivisions, which preserve the causal structure, as rigorously formalized by Dubut's pioneering work on natural systems (in Baues-Wirsching's sense \cite{CohomologySmallCategories}) associated with a directed space \cite{dubut_PhD,NaturalHomology}, and further developed in \cite{reversal-time-invariance}.

The globular analogue of this framework is the notion of cellular multipointed $d$-space (see Section~\ref{section:cellular-branching-space}), which consists of multipointed $d$-spaces constructed from globular cells, paired with the notion of globular subdivision (see Definition~\ref{def:globular-sbd}). This globular framework is sufficiently general to encompass all examples from concurrency theory, as precubical sets can be realized as cellular multipointed $d$-spaces (specifically by considering a cellular replacement of the realization of a precubical set as a multipointed $d$-space as defined e.g. in \cite[Definition~3.12]{GlobularNaturalSystem}). Following \cite{GlobularNaturalSystem}, globular analogues of Dubut's work confirm that globular subdivisions preserve the causal structure.

This paper serves as the first part of a two-part study focused on the explicit construction of the branching space and branching homology (along with the dual notions of merging space and merging homology after time-reversal). These homology theories encode information about the non-deterministic branching and merging areas of execution paths, thus capturing local aspects of the causal structure. Unlike previous approaches, this paper and its companion \cite{cubical-branching} aim to bypass the category of flows and the existing framework for branching homology introduced in \cite{exbranching,3eme}.

The first part, presented here, is dedicated to the globular setting, specifically focusing on cellular multipointed $d$-spaces and globular subdivisions. The second part \cite{cubical-branching} addresses the cubical analogue, exploring the setting of precubical sets and cubical subdivisions.

The primary goal of this paper is to provide an explicit geometric description of the branching space of a q-cofibrant multipointed $d$-space, achieved without relying on the categorization functor $\dcat:\ptop{\mathcal{M}}\to \dtop$ from multipointed $d$-spaces to flows. This construction serves as the globular analogue of the homotopy branching space of a precubical set, a concept fully developed in \cite{cubical-branching}, which provides a cubical perspective on branching spaces, offering a complementary approach to the globular setting explored here. The notion of a branching space naturally admits both globular and cubical formulations (much like other concepts in DAT). However, how to unify these two perspectives remains an open question.

As an illustrative application, we present a direct, purely topological (and indeed elementary) proof of the invariance of branching homology for cellular multipointed $d$-spaces under globular subdivision. While this result is not new (it can be derived from the main result of \cite{MultipointedSubdivision}, which establishes that two cellular multipointed $d$-spaces related by a zigzag of globular subdivisions have dihomotopy equivalent associated flows in the sense of \cite{hocont} via the categorization functor $\dcat:\ptop{\mathcal{M}} \to \dtop$, combined with \cite[Corollary~11.3]{3eme}), our proof offers a simpler and more intuitive approach. A similar invariance result is obtained in the cubical setting in \cite{cubical-branching}. However, this cubical analogue currently lacks alternative proofs, as the cubical analogue of the main result of \cite{MultipointedSubdivision} remains unproven.

The main results of this paper are summarized as follows. 

\begin{thmN} (Theorem~\ref{thm:dgerm-cgerm})
	Let $X$ be a q-cofibrant multipointed $d$-space. For all $\alpha\in X^0$, there is a homeomorphism $\mathcal{G}^-_\alpha(X) \iso \P^-_\alpha\dcat(X)$ between the branching space of $X$ at $\alpha$ introduced in Definition~\ref{def:branching-mdtop} and the branching space of the associated flow $\dcat(X)$ at $\alpha$ introduced in \cite{exbranching,3eme} and recalled in Definition~\ref{def:branching-flow}.
\end{thmN}

The letter $\mathcal{G}$ is used because the elements of $\mathcal{G}^-_\alpha(X)$ are germs of short directed paths starting from $\alpha$. The branching homology of a multipointed $d$-space $X$ is formally defined in Definition~\ref{def:hombrdef-mdtop} using the left derived functor of $\dcat:\ptop{\mathcal{M}} \to \dtop$. Theorem~\ref{thm:dgerm-cgerm} leads to a purely topological definition of the branching homology of a q-cofibrant multipointed $d$-space, without any use of the functor $\dcat:\ptop{\mathcal{M}}\to \dtop$, in Corollary~\ref{cor:gooddef-branching-homology-mdtop}.

\begin{thmN} (Theorem~\ref{thm:noncontracting3})
	Any globular subdivision \begin{tikzcd}[cramped]f:X\arrow[r,"\sbd"]&Y\end{tikzcd} induces a homeomorphism $\mathcal{G}^-_\alpha(X) \iso \mathcal{G}^-_{f(\alpha)}(Y)$ for all $\alpha\in X^0$ between the branching spaces of $X$ and $Y$ at $\alpha$.
\end{thmN}

\begin{thmN} (Corollary~\ref{cor:enfinfin})
	Let $f:X\to Y$ be a globular subdivision. Then for all $n\geq 0$, there are the isomorphisms $H_n^-(f):H_n^-(X)\longrightarrow H_n^-(Y)$ between the $n$-th branching homology groups.
\end{thmN}

The same definitions and results can be formulated for the merging space and the merging homology. Basically, it suffices to reverse the direction of time as follows. From a multipointed $d$-space $X$ (see Definition~\ref{def:multipointed-d-space}), we associate the opposite multipointed $d$-space $X^{op}$ defined as follows (with $c(t)=1-t$): the underlying topological space and the set of states are the same as $X$, and $\P_{\alpha,\beta}^{top}X^{op}=\{\gamma\in \ttop([0,1],|X|)\mid \gamma c\in \P_{\beta,\alpha}^{top}X\}$. The merging space of $X$ at $\alpha$ is defined by $\mathcal{G}^+_\alpha(X)=\mathcal{G}^-_\alpha(X^{op})$ and the merging homology $H^+_n(X)$ is equal to the branching homology $H^-_n(X^{op})$ for all $n\geq 0$. In particular, $H^+_0(X)$ is the free abelian group generated by the final states of $X^{op}$, which are the initial states of $X$. The details are left to the reader.

\subsection*{Outline of the paper}

Section~\ref{section:closed-incl} collects some facts about closed inclusions and $\Delta$-inclusions of $\Delta$-generated spaces. 

Section~\ref{section:branching-homology-mdtop} starts with some reminders about multipointed $d$-spaces, flows, the functor $\dcat$ from the former to the latter and their q-model structure. Definition~\ref{def:hombrdef-mdtop} introduces formally the branching homology of a multipointed $d$-space $X$ as the branching homology of the flow $\mathbf{L}\dcat(X)$ which is the image of $X$ under the left derived functor of $\dcat$.

Section~\ref{section:branching-mdtop} introduces the branching space of a multipointed $d$-space in Definition~\ref{def:branching-mdtop} after introducing the notion of short directed path in Definition~\ref{def:short}. The notion of short directed path is required to be able to use Corollary~\ref{cor:recognizing-closed-inclusions-3} in the proof of Proposition~\ref{prop:preparation}. 

Section~\ref{section:cellular-branching-space} studies the case of cellular multipointed $d$-spaces. The main results are Proposition~\ref{prop:G-colim} and Corollary~\ref{cor:push-glob} which describe the behavior of the branching space with respect to a cellular decomposition of a cellular multipointed $d$-space. As a first consequence, we prove in Corollary~\ref{cor:germ-glob} that the branching space of a q-cofibrant multipointed $d$-space is q-cofibrant, and in particular $\Delta$-Hausdorff.

Section~\ref{section:comparing-mdtop-flow} culminates with Theorem~\ref{thm:dgerm-cgerm} which establishes that the branching space of any q-cofibrant multipointed $d$-space is homeomorphic to the branching space of the associated flow.

Finally, Section~\ref{section:glb-sbd}, which contains the main application of this paper, establishes in Theorem~\ref{thm:noncontracting3} that globular subdivisions preserve the branching space of a cellular multipointed $d$-space using the topological material expounded in this paper. Note that we need to establish Theorem~\ref{thm:elm-sbd} which, roughly speaking, decomposes a globular subdivision in ``elementary'' globular subdivisions of the form $\globM(\mathbf{D}^{n})\to \globM(\mathbf{D}^{n})_F$. It could have been put in \cite{MultipointedSubdivision}: the latter theorem does not use the topological material of this paper.

\subsection*{Acknowledgments}

I am grateful to Tyrone Cutler for the observation in \cite{Tyrone}, which led to the proof of Proposition~\ref{prop:P-closedinclusion}.

\section{Closed inclusion}
\label{section:closed-incl}

In this paper, compact means Hausdorff quasi-compact. The category $\top$ denotes either the category of \textit{$\Delta$-generated spaces} or the category of \textit{$\Delta$-Hausdorff $\Delta$-generated spaces} (cf. \cite[Section~2 and Appendix~B]{leftproperflow}). It is Cartesian closed by a result due to Dugger and Vogt recalled in \cite[Proposition~2.5]{mdtop} and locally presentable by \cite[Corollary~3.7]{FR}. The internal hom is denoted by $\ttop(-,-)$. The right adjoint $k_\Delta:\mathcal{TOP}\to \top_\Delta$ of the inclusion from $\Delta$-generated spaces to general topological spaces is called the $\Delta$-kelleyfication. Let $\phi_U:k_\Delta(U)\to U$ be the counit. It preserves the underlying sets. One has $\ttop(-,-)=k_\Delta(\ttop_{co}(-,-))$ where $\ttop_{co}(-,-)$ means the set of continuous maps equipped with the compact-open topology. Every open subset of a $\Delta$-generated space equipped with the relative topology is $\Delta$-generated. A quotient map is a continuous map $f:X\to Y$ of $\Delta$-generated spaces which is onto and such that $Y$ is equipped with the final topology. The space $Y$ is called a final quotient of $X$. Every $\Delta$-generated space is sequential. Let $n\geq 1$. Denote by $\mathbf{D}^n = \{(x_1,\dots,x_n)\in \mathbb{R}^n, x_1^2+\dots + x_n^2 \leq 1\}$ the $n$-dimensional disk, and by $\mathbf{S}^{n-1} = \{(x_1,\dots,x_n)\in \mathbb{R}^n, x_1^2+\dots + x_n^2 = 1\}$ the $(n-1)$-dimensional sphere. By convention, let $\mathbf{D}^{0}=\{0\}$ and $\mathbf{S}^{-1}=\varnothing$.

\bd \label{def:closed-incl} A \textit{closed inclusion} $f:A\to B$ of $\Delta$-generated spaces is a one-to-one continuous map such that $f(A)$ is a closed subset of $B$ and such that $f$ induces a homeomorphism between $A$ and $f(A)$ equipped with the relative topology. 
\ed

Every $h$-cofibration of $\top$ is a closed inclusion. This implies that $f(A)$ equipped with the relative topology is $\Delta$-generated. Recall that a closed subset of a $\Delta$-generated space equipped with the relative topology is not necessarily $\Delta$-generated. However, it is always sequential, being a closed subset of a sequential space. For example, the Cantor set is a closed subset of the $\Delta$-generated space $[0,1]$ which is not $\Delta$-generated for the relative topology, the latter having a discrete $\Delta$-Kelleyfication since it is totally disconnected. 

\bd \label{def:delta-incl}
A one-to-one map of $\Delta$-generated spaces $f:A\to B$ is a \textit{$\Delta$-inclusion} if for all $\Delta$-generated spaces $Z$, the set map $Z\to A$ is continuous if and only if the composite set map $Z\to A\to B$ is continuous.
\ed

\begin{prop} \label{DeltaIncl} Let $f:A\to B$ be a one-to-one continuous map. The following assertions are equivalent:
\begin{enumerate}
	\item $f$ is a $\Delta$-inclusion.
	\item $A$ is homeomorphic to $f(A)$ equipped with the $\Delta$-kelleyfication of the relative topology.
	\item A set map $[0,1]\to A$ is continuous if and only if the composite set map $[0,1]\to A\to B$ is continuous.
\end{enumerate}
\end{prop}

\bpf The proof is similar to the same statement for $k$-inclusions of $k$-spaces.
\epf

Let $A$ be an open subset of a $\Delta$-generated space $B$. Then the inclusion $A\subset B$ is a $\Delta$-inclusion, $A$ being equipped with the relative topology. Every closed inclusion of $\Delta$-generated spaces is a $\Delta$-inclusion with closed image. The converse is false. The one-to-one map $\{1/n\mid n\geq 1\} \cup \{0\} \subset [0,1]$ with the left-hand space equipped with the discrete topology is a $\Delta$-inclusion with closed image. However, the relative topology of $\{1/n\mid n\geq 1\}\cup \{0\}$ is not discrete.

\bth \label{thm:recognizing-closed-inclusions-1}
Let $f:A\to B$ be a one-to-one continuous map between $\Delta$-generated spaces. Suppose that $B$ is $\Delta$-Hausdorff. The map $f$ is a closed inclusion if and only if for any sequence $(x_n)_{n\geq 0}$ of $A$, if the sequence $(f(x_n))_{n\geq 0}$ of $B$ is convergent, so is the sequence $(x_n)_{n\geq 0}$ of $A$.
\eth

\bpf Since $B$ is $\Delta$-Hausdorff, so is $A$. Thus, both $A$ and $B$ have unique sequential limits by \cite[Proposition~B.17]{leftproperflow}. Suppose that $f$ is a closed inclusion. It means that there is a homeomorphism $A\iso f(A)$ between $A$ and $f(A)$ equipped with the relative topology. Let $(x_n)_{n\geq 0}$ be a sequence of $A$ such that the sequence $(f(x_n))_{n\geq 0}$ of $B$ is convergent. Since $f(A)$ is closed in $B$, the limit is of the form $f(x_\infty)$. Because of the homeomorphism $A\iso f(A)$, we deduce that $x_\infty$ is the limit of $(x_n)_{n\geq 0}$. Conversely, let $f:A\to B$ be a map satisfying the condition of the theorem. Let $y_\infty$ be an element of the closure $\cl(f(A))$ of $f(A)$ in $B$. Then $y_\infty$ is the limit of some sequence $(f(x_n))_{n\geq 0}$ of $f(A)$ because $B$ is sequential, being a $\Delta$-generated space. We deduce that the sequence $(x_n)_{n\geq 0}$ of $A$ is convergent to, say, $x_\infty$. By continuity of $f$ and by uniqueness of the sequential limits, we obtain $y_\infty=f(x_\infty)$. We have proved that $\cl(f(A))\subset  f(A)$. We deduce that $\cl(f(A))=  f(A)$, i.e. that $f(A)$ is a closed subset of $B$. Consequently, equipped with the relative topology, $f(A)$ is sequential, $B$ being sequential. The space $A$ is sequential, being $\Delta$-generated. The condition of the theorem means that the sequential spaces $A$ and $f(A)$ have the same convergent sequences. This means that $f$ induces a homeomorphism $A \iso f(A)$. We have proved that $f$ is a closed inclusion.
\epf

The hypothesis $\Delta$-Hausdorff is necessary. The map $\{0\}\subset \{0^-,0^+\}$ where $\{0^-,0^+\}$ is equipped with the indiscrete topology satisfies the condition of Theorem~\ref{thm:recognizing-closed-inclusions-1} and it is not a closed inclusion.

Proposition~\ref{prop:quotient-closedinclusion} is a standard fact in algebraic topology textbooks for $k$-spaces.

\begin{prop} \label{prop:quotient-closedinclusion} (\cite[Proposition~B.14]{leftproperflow}) Consider the commutative diagram of $\Delta$-generated spaces
\[
\begin{tikzcd}[column sep=3em,row sep=3em]
	X \arrow[r,"p"] \arrow[d,"f"'] & Z \arrow[d,"\overline{f}"] \\
	Y \arrow[r,"q"] & W
\end{tikzcd}
\]
If $f$ is a closed inclusion, $\overline{f}$ is one-to-one, $p$ is onto and either $q$ is closed (i.e. the image of a closed subset is a closed subset) or $q$ is a quotient map, with $q^{-1}(\overline{f}(Z))\subset f(X)$, then $\overline{f}$ is a closed inclusion.
\end{prop}

\begin{cor} \label{cor:recognizing-closed-inclusions-3}
	Consider a closed inclusion $f:X\to Y$ of $\Delta$-generated spaces. Assume that $\mathcal{R}$ is an equivalence relation on $X$ and $\mathcal{S}$ an equivalence relation on $Y$ such that $x_1\mathcal{R}x_2$ if and only if $f(x_1)\mathcal{S}f(x_2)$. Equip the quotient sets $X/\mathcal{R}$ and $Y/\mathcal{S}$ with the final topology. Then $\overline{f}$ induces a continuous one-to-one map $\overline{f}:X/\mathcal{R}\to Y/\mathcal{S}$ and there is a commutative diagram of spaces 
	\[
	\begin{tikzcd}[column sep=3em,row sep=3em]
		X \arrow[r,"p"] \arrow[d,"f"'] & X/\mathcal{R} \arrow[d,"\overline{f}"] \\
		Y \arrow[r,"q"] & Y/\mathcal{S}
	\end{tikzcd}
	\]
	Assume that $q^{-1}(\overline{f}(X/\mathcal{R}))\subset f(X)$. If $f$ is a closed inclusion, then $\overline{f}$ is a closed inclusion as well.
\end{cor}

\bpf From the hypothesis $x_1\mathcal{R}x_2$ if and only if $f(x_1)\mathcal{S}f(x_2)$, we deduce that $\overline{f}$ is one-to-one because the quotient spaces $X/\mathcal{R}$ and $Y/\mathcal{S}$ are equipped with the final topology which does not change the underlying sets. The proof is complete thanks to Proposition~\ref{prop:quotient-closedinclusion}.
\epf

\begin{cor} \label{cor:recognizing-closed-inclusions-4}
	Consider a closed inclusion $f:X\to Y$ of $\Delta$-generated spaces. Assume that $\mathcal{S}$ is an equivalence relation on $Y$. Equip the quotient set $Y/\mathcal{S}$ with the final topology. Assume that there is a commutative diagram of spaces
	\[
	\begin{tikzcd}[column sep=3em,row sep=3em]
		X \arrow[r,"p"] \arrow[d,"f"'] & Z \arrow[d,"\overline{f}"] \\
		Y \arrow[r,"q"] & Y/\mathcal{S}
	\end{tikzcd}
	\]
	such that $p$ is onto, $f$ is a closed inclusion and $\overline{f}$ is a continuous bijection. Assume that $q^{-1}(\overline{f}(Z))\subset f(X)$. Then $\overline{f}$ is a homeomorphism.
\end{cor}

\bpf
By Proposition~\ref{prop:quotient-closedinclusion}, $\overline{f}$ is a closed inclusion. Therefore it is a homeomorphism, being a continuous bijection by hypothesis.
\epf

\section{Branching homology of a multipointed d-space}
\label{section:branching-homology-mdtop}

Let $U$ be a topological space. A \textit{(Moore) path} of $U$ consists of a continuous map $\gamma:[0,\ell]\to U$. The real number $\ell$ is called the \textit{length} of $\gamma$. Let $\gamma_1:[0,\ell_1]\to U$ and $\gamma_2:[0,\ell_2]\to U$ be two Moore paths of a topological space $U$ such that $\gamma_1(\ell_1)=\gamma_2(0)$. The \textit{Moore composition} $\gamma_1*\gamma_2:[0,\ell_1+\ell_2]\to U$ is the Moore path defined by 
\[
(\gamma_1*\gamma_2)(t)=
\begin{cases}
	\gamma_1(t) & \hbox{ for } t\in [0,\ell_1]\\
	\gamma_2(t-\ell_1) &\hbox{ for }t\in [\ell_1,\ell_1+\ell_2].
\end{cases}
\]
The Moore composition of Moore paths is strictly associative. Let $\gamma_1$ and $\gamma_2$ be two continuous maps from $[0,1]$ to some topological space such that $\gamma_1(1)=\gamma_2(0)$. The continuous map defined by 
\[
(\gamma_1 *_N \gamma_2)(t) = \begin{cases}
	\gamma_1(2t)& \hbox{ if }0\leq t\leq \frac{1}{2},\\
	\gamma_2(2t-1)& \hbox{ if }\frac{1}{2}\leq t\leq 1
\end{cases}
\]
is called the \textit{normalized composition}. 

The set $\mathcal{M}$ of non-decreasing surjective maps from $[0,1]$ to $[0,1]$ is equipped with the $\Delta$-kelleyfication of the relative topology induced by the inclusion $\mathcal{M} \subset \ttop([0,1],[0,1])$. 

\bd \cite[Definition~3.4]{Moore3} \label{def:multipointed-d-space} A \textit{multipointed $d$-space} $X$ is a triple $(|X|,X^0,\P^{top}X)$ such that
\begin{itemize}
	\item The pair $(|X|,X^0)$ is a \textit{multipointed space}. The space $|X|$ is called the \textit{underlying space} of $X$ and the set $X^0$ the \textit{set of states} of $X$.
	\item The set $\P^{top}X$ is a set of continuous maps from $[0,1]$ to $|X|$ called the \textit{execution paths}, satisfying the following axioms:
	\begin{itemize}
		\item For any execution path $\gamma$, one has $\gamma(0),\gamma(1)\in X^0$.
		\item Let $\gamma$ be an execution path of $X$. Then any composite $\gamma\phi$ with $\phi\in \mathcal{M}$ is an execution path of $X$.
		\item Let $\gamma_1$ and $\gamma_2$ be two execution paths of $X$; if the normalized composition $\gamma_1 *_N \gamma_2$ exists, it is an execution path of $X$.
	\end{itemize}
\end{itemize}
A map $f:X\to Y$ of multipointed $d$-spaces is a map of multipointed spaces from $(|X|,X^0)$ to $(|Y|,Y^0)$ such that for any execution path $\gamma$ of $X$, the map $\P^{top}f:\gamma\mapsto f. \gamma$ is an execution path of $Y$. The category of multipointed $d$-spaces is denoted by $\ptop{\mathcal{M}}$. Let $\P_{\alpha,\beta}^{top} X = \{\gamma\in \P^{top}X\mid \gamma(0)=\alpha,\gamma(1)=\beta\}$. The set $\P^{top}_{\alpha,\beta} X$ is equipped with the $\Delta$-kelleyfication of the relative topology with respect to the inclusion $\P^{top}_{\alpha,\beta} X \subset \ttop([0,1],|X|)$. Thus a set map $Z\to \P^{top}_{\alpha,\beta} X$ where $Z$ is $\Delta$-generated, is continuous if and only if the associated map $Z\p [0,1] \to |X|$ is continuous. 
\ed

The category $\ptop{\mathcal{M}}$ is locally presentable by \cite[Proposition~6]{Moore3}.  Every set $S$ can be viewed as a multipointed $d$-space $(S,S,\varnothing)$. The \textit{topological globe of a topological space $Z$}, which is denoted by $\globM(Z)$, is the multipointed $d$-space defined as follows
\begin{itemize}
	\item the underlying topological space is the quotient space \[\frac{\{{0},{1}\}\sqcup (Z\p[0,1])}{(z,0)=(z',0)={0},(z,1)=(z',1)={1}}\]
	\item the set of states is $\{{0},{1}\}$
	\item the set of execution paths is the set of continuous maps \[\{\delta_z\phi\mid \phi\in \mathcal{M},z\in  Z\}\]
	with $\delta_z(t) = (z,t)$.	It is equal to the underlying set of the space $Z \p \mathcal{M}$.
\end{itemize}
In particular, $\globM(\varnothing)$ is the multipointed $d$-space $\{{0},{1}\} = (\{{0},{1}\},\{{0},{1}\},\varnothing)$. Let $\vI^{top} = \globM(\mathbf{D}^0)$.

The \textit{q-model structure} of multipointed $d$-spaces \cite[Section~4]{Moore3} is the unique combinatorial model structure such that $I^{gl,top}\cup \{C:\varnothing \to \{0\},R:\{0,1\} \to \{0\}\}$ with $I^{gl,top}=\{\globM(\mathbf{S}^{n-1})\subset \globM(\mathbf{D}^{n}) \mid n\geq 0\}$ is the set of generating cofibrations, the maps between globes being induced by the closed inclusions $\mathbf{S}^{n-1}\subset \mathbf{D}^{n}$, and such that $J^{gl,top}=\{\globM(\mathbf{D}^{n})\subset \globM(\mathbf{D}^{n+1}) \mid n\geq 0\}$ is the set of generating trivial cofibrations, the maps between globes being induced by the closed inclusions $(x_1,\dots,x_n)\mapsto (x_1,\dots,x_n,0)$. The weak equivalences are the maps of multipointed $d$-spaces $f:X\to Y$  inducing a bijection $f^0:X^0\iso Y^0$ and a weak homotopy equivalence $\P^{{top}} f:\P^{{top}}_{\alpha,\beta} X \to \P^{{top}}_{f(\alpha),f(\beta)} Y$ for all $(\alpha,\beta)\in X^0\p X^0$ and the fibrations are the maps of multipointed $d$-spaces $f:X\to Y$  inducing a q-fibration $\P_{\alpha,\beta}^{{top}} f:\P_{\alpha,\beta}^{{top}} X \to \P_{f(\alpha),f(\beta)}^{{top}} Y$ of topological spaces for all $(\alpha,\beta)\in X^0\p X^0$. All multipointed $d$-spaces are q-fibrant.

\bd \cite[Definition~4.11]{model3} \label{def:flow}
A \textit{flow} $X$ is a small topologically enriched semicategory. Its set of objects (preferably called \textit{states}) is denoted by $X^0$ and the space of morphisms (preferably called \textit{execution paths}) from $\alpha$ to $\beta$ is denoted by $\P_{\alpha,\beta}X$. For any $x\in \P_{\alpha,\beta}X$, let $s(x)=\alpha$ and $t(x)=\beta$. The category is denoted by $\dtop$. 
\ed

The category $\dtop$ is locally presentable. Every set can be viewed as a flow with an empty path space. This gives rise to a functor from sets to flows which is limit-preserving and colimit-preserving. More generally, any poset can be viewed as a flow, with a unique execution path from $u$ to $v$ if and only if $u<v$. This gives rise to a functor from the category of posets together with the strictly increasing maps to flows.

\begin{nota} \label{nota:glob}
	For any topological space $Z$, the flow $\glob(Z)$ is the flow having two states $0$ and $1$ and such that the only nonempty space of execution paths, when $Z$ is nonempty, is $\P_{0,1}\glob(Z)=Z$. It is called \textit{the globe of $Z$}. Let $\vI = \glob(\mathbf{D}^0)$.
\end{nota}

The \textit{q-model structure} of flows \cite[Theorem~7.6]{QHMmodel} is the unique combinatorial model structure such that 
\[I^{gl} \cup \{C:\varnothing \to \{0\},R:\{0,1\} \to \{0\}\}\]
with $I^{gl}=\{\glob(\mathbf{S}^{n-1})\subset \glob(\mathbf{D}^{n}) \mid n\geq 0\}$ is the set of generating cofibrations, the maps between globes being induced by the closed inclusions $\mathbf{S}^{n-1}\subset \mathbf{D}^{n}$, and such that 
\[
J^{gl}=\{\glob(\mathbf{D}^{n})\subset \glob(\mathbf{D}^{n+1}) \mid n\geq 0\}
\]
is the set of generating trivial cofibrations, the maps between globes being induced by the closed inclusions $(x_1,\dots,x_n)\mapsto (x_1,\dots,x_n,0)$. The weak equivalences are the maps of flows $f:X\to Y$  inducing a bijection $f^0:X^0\iso Y^0$ and a weak homotopy equivalence  $\P f:\P_{\alpha,\beta} X \to \P_{f(\alpha),f(\beta)} Y$ for all $(\alpha,\beta)\in X^0\p X^0$ and the fibrations are the maps of flows $f:X\to Y$  inducing a q-fibration $\P_{\alpha,\beta} f:\P_{\alpha,\beta} X \to \P_{f(\alpha),f(\beta)} Y$ of topological spaces for all $(\alpha,\beta)\in X^0\p X^0$. All flows are q-fibrant.

\bth \label{thm:def-Cminus} \cite[Theorem~5.5]{exbranching}
Let $Z$ be a topological space. The data $C^-(Z)^0=\{0\}$, $\P_{0,0}Z=Z$ and $x*y=x$ assemble to a flow $C^-(Z)$. The mapping $Z\mapsto C^-(Z)$ gives rise to a right Quillen adjoint $C^-:\top\to \dtop$ from the q-model structure of $\top$ to the q-model structure of $\dtop$. 
\eth

The left Quillen adjoint $\P^-:\dtop \to \top$ is defined by taking a flow $X$ to the quotient space 
\[
\P^- X = \coprod_{\alpha\in X^0} \P^-_\alpha X
\]
where $\P_\alpha X$ is the quotient of the space $\coprod_{\beta\in X^0} \P_{\alpha,\beta}X$ by the equivalence relation generated by the identifications $u*v=u$. By \cite[Theorem~4.1]{exbranching}, the branching space functor $\P^-:\dtop \to \top$ is badly behaved with respect to the weak equivalences. Indeed, there exists a weak equivalence of flows $f:X\to Y$ such that $\P^-f:\P^-X\to \P^-Y$ is not a weak homotopy equivalence of spaces. Hence the following definition.

\bd \label{def:branching-flow}
Let $X$ be a flow. The topological space $\P^-X$ is called the \textit{branching space} of $X$. The topological space $\hop^-X=\P^- X^{cof}$ is called the \textit{homotopy branching space} of $X$ where $(-)^{cof}$ is some q-cofibrant replacement of $X$. The latter space is unique only up to homotopy equivalence.
\ed

There is a unique functor $\dcat:\ptop{\mathcal{M}}\longrightarrow \dtop$ from the category of multipointed $d$-spaces to the category of flows, called the \textit{categorization functor}, taking a multipointed $d$-space $X$ to the unique flow $\dcat(X)$ such that $\dcat(X)^0=X^0$ and such that $\P_{\alpha,\beta}\dcat(X)$ is the quotient of the space of execution paths $\P_{\alpha,\beta}^{top}X$ by the equivalence relation generated by the reparametrization by $\mathcal{M}$, the composition of $\dcat(X)$ being induced by the normalized composition. One has $\dcat(\vI^{top})=\vI$. For any topological space $Z$, there is the natural isomorphism of flows $\dcat(\globM(Z))\iso \glob(Z)$.

We recall the definition of the branching homology of flows. 

\bd\label{def:hombrdef-flow}  \cite[Definition~6.1]{exbranching}
Let $X$ be a flow. The \textit{branching homology} groups $H_{*}^-(X)$ are defined as follows: 
\begin{enumerate}
	\item for $n\geq 1$, $H_{n+1}^-(X):=H_n(\hop^-(X))$
	\item  $H_1^-(X):=\ker(\epsilon)/\im(\partial)$
	\item $H_0^-(X):=\mathbb{Z}[X^0]/\im(\epsilon)$
\end{enumerate}
with the augmentation $\epsilon:\mathbb{Z}[\hop^-(X)]\to \mathbb{Z}[X^0]$ defined by $\epsilon(\gamma)=\gamma(0)$ and $\partial:\mathbb{Z}[\top([0,1],\hop^-(X))]\to \mathbb{Z}[\hop^-(X)]$ defined by $\partial(f)=f(0)-f(1)$, where $\ker$ is the kernel and $\im$ the image. For any flow $X$, $H_0^-(X)$ is the free abelian group generated by the final states of $X$. The branching homology of $X$ is clearly $X^0$-graded. Let 
\[
H^-_*(X)= \coprod_{\alpha\in X^0} G^\alpha H^-_*(X).
\]
\ed

By \cite[Corollary~6.5]{exbranching}, for every weak equivalence of flows $f:X\to Y$, the morphism of abelian groups $H_n^-(f):H_n^-(X)\to H_n^-(Y)$ is an isomorphism for all $n\geq 0$. By \cite[Corollary~11.3]{3eme}, for any retract of a transfinite composition of pushouts of generating subdivision in the sense of \cite[Definition~8.3]{MultipointedSubdivision} $f:X\to Y$, the morphism of abelian groups $H_n^-(f):H_n^-(X)\to H_n^-(Y)$ is an isomorphism for all $n\geq 0$.

The functor $\dcat:\ptop{\mathcal{M}}\longrightarrow \dtop$ is neither a left nor a right adjoint. However, by \cite[Theorem~15]{Moore3}, it takes q-cofibrant multipointed $d$-spaces to q-cofibrant flows and its left derived functor $\mathbf{L}\dcat=\dcat\circ(-)^{cof}$ (where $(-)^{cof}$ is a q-cofibrant replacement functor of the q-model structure of multipointed $d$-spaces) in the sense of \cite{HomotopicalCategory} induces an equivalence of categories between the homotopy categories of the q-model structures of multipointed $d$-spaces and flows. Moreover, the functor $\dcat:\ptop{\mathcal{M}}\longrightarrow \dtop$ takes weak equivalences between q-cofibrant multipointed $d$-spaces to weak equivalences between q-cofibrant flows. We can now define the branching homology of a multipointed $d$-space as follows.

\bd\label{def:hombrdef-mdtop}
Let $X$ be a multipointed $d$-space. Let $n\geq 0$. The \textit{$n$-th branching homology group} $H_n^-(X)$ is the branching homology group $H_n^-(\mathbf{L}\dcat(X))$. 
\ed 

Definition~\ref{def:hombrdef-mdtop} makes sense since the flow $\mathbf{L}\dcat(X)$ is unique up to weak equivalences. By \cite[Theorem~9.3]{MultipointedSubdivision}, any globular subdivision $f:X\to Y$ in the sense of Definition~\ref{def:globular-sbd}, $\mathbf{L}\dcat(f)=\dcat(f):\dcat(X)\to \dcat(Y)$ is the composite of a transfinite composition of pushouts of generating subdivision in the sense of \cite[Definition~8.3]{MultipointedSubdivision} and of a weak equivalence of flows. Therefore, by \cite[Corollary~11.3]{3eme} and \cite[Corollary~6.5]{exbranching} mentioned above, the induced map of abelian groups $H_n^-(f):H_n^-(X)\to H_n^-(Y)$ is an isomorphism for all $n\geq 0$.

\section{Branching space of a multipointed d-space}
\label{section:branching-mdtop}

The set $\mathcal{I}$ of non-decreasing continuous maps from $[0,1]$ to $[0,1]$ is equipped with the $\Delta$-kelleyfication of the relative topology induced by the inclusion $\mathcal{I} \subset \ttop([0,1],[0,1])$.

\bd \cite[Definition~1.1]{mg} \cite[Definition~4.1]{DAT_book} \label{def:directed_space}
A \textit{directed space} is a pair $X=(|X|,\vec{P}(X))$ consisting of a topological space $|X|$ and a set $\vec{P}(X)$ of continuous paths from $[0,1]$ to $|X|$ satisfying the following axioms:
\begin{itemize}
	\item $\vec{P}(X)$ contains all constant paths;
	\item $\vec{P}(X)$ is closed under normalized composition;
	\item $\vec{P}(X)$ is closed under reparametrization by an element of $\mathcal{I}$.
\end{itemize}
The space $|X|$ is called the \textit{underlying topological space} or the \textit{state space}. The elements of $\vec{P}(X)$ are called \textit{directed paths}. A morphism of directed spaces is a continuous map between the underlying topological spaces which takes a directed path of the source to a directed path of the target. The category of directed spaces is denoted by $\ptop{}$. The space $\vec{P}(X)$ is equipped with the $\Delta$-kelleyfication of the compact-open topology.
\ed

By \cite[Proposition~3.7 and Theorem~3.8]{GlobularNaturalSystem}, the mapping $\discont:Y=(|Y|,\vec{P}(Y)) \mapsto (|Y|,|Y|,\vec{P}(Y))$ induces a full and faithful functor $\discont:\ptop{}\to\ptop{\mathcal{M}}$ which is a right adjoint. By \cite[Proposition~3.6]{GlobularNaturalSystem}, the left adjoint $\cont:\ptop{\mathcal{M}} \to \ptop{}$ is defined as follows. The underlying space of $\cont(X)$ is $|X|$ and the set of directed paths $\vec{P}(\cont(X))$ consists of all constant paths and all Moore compositions of the form $(\gamma_1\phi_1\mu_{\ell_1}) * \dots * (\gamma_n\phi_n\mu_{\ell_n})$ such that $\ell_1+\dots + \ell_n = 1$ where $\gamma_1,\dots,\gamma_n$ are execution paths of $X$ and $\phi_i\in \mathcal{I}$ for $i=1,\dots,n$, and where $\mu_\ell:[0,\ell]\to [0,1]$ is defined by $\mu_\ell(t)=t/\ell$ with $\ell>0$.

\bd \label{def:short}
Let $X$ be a multipointed $d$-space. A \textit{directed path} of $X$ is a directed path of $\cont(X)$. A \textit{short} directed path of $X$ is a directed path $\gamma:[0,1]\to |X|$ of $X$ such that $\gamma(]0,1[)\cap X^0=\varnothing$. The set of short directed paths of $X$ starting from $\alpha$ is denoted by $\mathcal{P}^-_\alpha(X)$. It is equipped with the $\Delta$-kelleyfication of the relative topology for the inclusion $\mathcal{P}^-_\alpha(X) \subset \ttop([0,1],|X|)$. Let 
\[
\mathcal{P}^-(X) = \coprod_{\alpha\in X^0} \mathcal{P}^-_\alpha(X).
\]
\ed

\begin{nota}  \label{nota:Iminus-def}
	Let $\mathcal{I}^-\subset \mathcal{I}$ be the subspace of maps $\phi\in \mathcal{I}$ such that $\phi(0)=0$ and $\phi(]0,1[)\cap \{0,1\}=\varnothing$ equipped with the $\Delta$-kelleyfication of the relative topology.
\end{nota}

\bd  \label{def:branching-mdtop}
Let $X$ be a multipointed $d$-space. Let $\alpha\in X^0$. The \textit{branching space} of $X$ at $\alpha$ is the \textit{final} quotient 
\[
\mathcal{G}^-_\alpha(X) = \mathcal{P}^-_\alpha(X)/\!\simeq^- 
\]
where the equivalence relation $\simeq^-$ is generated by the identifications $\gamma_1= \gamma_2$ if and only if there exist $\gamma\in \mathcal{P}^-_\alpha(X)$ and $\phi_1,\phi_2\in \mathcal{I}^-$ such that $\gamma\phi_1=\gamma_1$ and $\gamma\phi_2=\gamma_2$. Let 
\[
\mathcal{G}^-(X) = \coprod_{\alpha\in X^0} \mathcal{G}^-_\alpha(X).
\]
The image of $\gamma$ under the canonical map $\mathcal{P}^-_\alpha(X)\to \mathcal{G}^-_\alpha(X)$ is denoted by $\tr{\gamma}^-$.
\ed

The space $\mathcal{G}^-_\alpha(X)$ is $\Delta$-Hausdorff if and only if the graph of the equivalence relation $\simeq^-$ is closed by \cite[Corollary~B.11]{leftproperflow}. By Corollary~\ref{cor:germ-glob}, this condition holds when $X$ is q-cofibrant. We have to work with short directed paths to be able to use Corollary~\ref{cor:recognizing-closed-inclusions-3} in the proof of Proposition~\ref{prop:preparation}. Because of this, the mapping $X\mapsto \mathcal{G}^-(X)$ is not functorial. Hence the importance of Proposition~\ref{prop:noncontracting1} and Corollary~\ref{cor:functoriality}.

\begin{prop} \label{prop:noncontracting1}
Let $n\geq 0$. Consider a pushout diagram of cellular multipointed $d$-spaces
\[
\begin{tikzcd}[row sep=3em,column sep=3em]
	\globM(\mathbf{S}^{n-1}) \arrow[r,"g"] \arrow[d,hookrightarrow] & A \arrow[d,hookrightarrow,"f"] \\
	\globM(\mathbf{D}^{n}) \arrow[r,"\widehat{g}"] & \cocartesian X
\end{tikzcd}
\]
Then the mapping $\gamma\mapsto f\gamma$ induces continuous maps $\mathcal{P}^-_\alpha(A) \to \mathcal{P}^-_\alpha(X)$ and $\mathcal{G}^-_\alpha(A) \to \mathcal{G}^-_\alpha(X)$ for all $\alpha\in A^0$. 
\end{prop}

\bpf
Let $\gamma\in \mathcal{P}^-_\alpha(A)$. The functor $(-)^0$ being colimit preserving, $f$ induces a bijection from $A^0$ to $X^0$. The underlying space functor $|-|:\ptop{\mathcal{M}}\to \top$ being a left adjoint by \cite[Proposition~A.2]{MultipointedSubdivision}, it is colimit preserving. The endofunctor $Z\mapsto |\globM(Z)|$ being a left Quillen adjoint for the q-model structure of $\top$ by \cite[Proposition~A.2]{MultipointedSubdivision}, the map $|f|:|A|\to |X|$ is one-to-one, being a q-cofibration. Thus $f(\gamma(]0,1]))\cap X^0=\varnothing$ and the proof is complete.
\epf

\begin{prop} \label{prop:globe-short-directed-path}
	Let $Z$ be a topological space. There are the homeomorphisms
	\[
	\mathcal{P}^-_\alpha(\globM(Z)) \iso  \begin{cases}
		Z \p \mathcal{I}^- & \hbox{ if }\alpha=0\\
		\varnothing& \hbox{ if }\alpha=1\\
	\end{cases}
	\]
\end{prop}

\bpf One has $\mathcal{P}^-_1(\globM(Z))=\varnothing$ because there is no nonconstant directed path in $\cont(\globM(Z))$ starting from $1$. The set map  
\[
\begin{tikzcd}[row sep=0em,column sep=3em]
	Z\p \mathcal{I}^- \arrow[r,"\Psi_1"] & \mathcal{P}^-_0(\globM(Z))\\
	(z,\phi) \arrow[r,mapsto] & \delta_z\phi
\end{tikzcd}
\]
is continuous since the mapping $(t,z,\phi)\mapsto (z,\phi(t))$ from $[0,1]\p Z\p \mathcal{I}^-$ to $|\globM(Z)|$ is continuous. Let $\pr_1:|\globM(Z)|\backslash\{0,1\}\to Z$ and $\pr_2:|\globM(Z)|\to |\globM(\{0\})|=[0,1]$ be the two projection maps. The set map 
\[
\begin{tikzcd}[row sep=0em,column sep=3em]
	\mathcal{P}^-_0(\globM(Z)) \arrow[r,"\Psi_2"] & Z\p \mathcal{I}^-  \\
	\gamma \arrow[r,mapsto] & (\pr_1\gamma(1/2),\pr_2\gamma)
\end{tikzcd}
\]
is continuous since the mapping $(t,\gamma)\mapsto (\pr_1\gamma(1/2),\pr_2\gamma(t))$ from $[0,1]\p \mathcal{P}^-_0(\globM(Z))$ to $Z\p [0,1]$ is continuous. Besides, $\Psi_1\Psi_2=\id$ and $\Psi_2\Psi_1=\id$. Hence the proof is complete.
\epf

\begin{cor} \label{cor:globe-branching}
	Let $Z$ be a topological space. There are the homeomorphisms
	\[
	\mathcal{G}^-_\alpha(\globM(Z)) \iso  \begin{cases}
		Z & \hbox{ if }\alpha=0\\
		\varnothing& \hbox{ if }\alpha=1\\
	\end{cases}
	\]
\end{cor}

\bpf One has $\mathcal{G}^-_1(\globM(Z))=\varnothing$ because there is no non-constant directed path in $\cont(\globM(Z))$ starting from $1$. For any $\phi_1,\phi_2\in \mathcal{I}^-$, one has $\id_{[0,1]}\phi_1=\phi_1$ and $\id_{[0,1]}\phi_2=\phi_2$. Hence the proof is complete using Proposition~\ref{prop:globe-short-directed-path}. 
\epf

\begin{cor} \label{cor:functoriality}
	Let $f:U\to V$ be a continuous map of topological spaces. The mapping $\gamma\mapsto \globM(f)\gamma$ induces continuous maps $\mathcal{P}^-_\alpha(\globM(U)) \to \mathcal{P}^-_\alpha(\globM(V))$ and $\mathcal{G}^-_\alpha(\globM(U)) \to \mathcal{G}^-_\alpha(\globM(V))$ for $\alpha\in \{0,1\}$.
\end{cor}

\section{The case of q-cofibrant multipointed d-spaces}
\label{section:cellular-branching-space}

For a class of maps $\C$, $\cell(\C)$ denotes the class of transfinite compositions of pushouts of elements of $\C$. A \textit{cellular} object $X$ of a cofibrantly generated model category is an object such that the canonical map $\varnothing\to X$ belongs to $\cell(I)$ where $I$ is the set of generating cofibrations. The maps of $\cell(I)$ are called \textit{cellular maps}. The transfinite sequence of pushouts is called a \textit{cellular decomposition} of $X$ and each pushout is called a \textit{cell}. A \textit{transfinite tower} in a cocomplete category $\K$ is a colimit-preserving functor $\lambda\to \K$ from an ordinal $\lambda$ viewed as a small category thanks to its poset structure to $\K$.

All cellular multipointed $d$-spaces for the q-model structure of $\ptop{\mathcal{M}}$ can be reached from $\varnothing$ without using the cofibration $R:\{0,1\}\to \{0\}$ and by regrouping the pushouts of $C:\varnothing\to \{0\}$ at the very beginning by \cite[Proposition~3.8]{MultipointedSubdivision}. Thus, for the sequel, a cellular decomposition of a cellular multipointed $d$-space $X$ of the q-model category $\ptop{\mathcal{M}}$ consists of a colimit-preserving functor $\widetilde{X}:\lambda \longrightarrow \ptop{\mathcal{M}}$ from a transfinite ordinal $\lambda$ to the category of multipointed $d$-spaces such that
\begin{itemize}
	\item The multipointed $d$-space $\widetilde{X}_0$ is a set, in other words $\widetilde{X}_0=(X^0,X^0,\varnothing)$ for some set $X^0$.
	\item For all $\nu<\lambda$, there is a pushout diagram of multipointed $d$-spaces 
	\[\begin{tikzcd}[row sep=3em, column sep=3em]
		\globM(\mathbf{S}^{n_\nu-1}) \arrow[d,rightarrowtail] \arrow[r,"g_\nu"] & \widetilde{X}_\nu \arrow[d,rightarrowtail] \\
		\globM(\mathbf{D}^{n_\nu}) \arrow[r,"\widehat{g_\nu}"] & \cocartesian \widetilde{X}_{\nu+1}
	\end{tikzcd}\]
	with $n_\nu \geq 0$
	\item $X=\liminj \widetilde{X}$.
\end{itemize}
The underlying topological space $|X|$ is Hausdorff by \cite[Proposition~4.4]{Moore3}. For all $\nu\leq \lambda$, there is the equality $\widetilde{X}_\nu^0=X^0$. Denote by $c_\nu = |\globM(\mathbf{D}^{n_\nu})|\backslash |\globM(\mathbf{S}^{n_\nu-1})|$ the $\nu$-th cell of $X$. It is called a \textit{globular cell}. Like in the usual setting of CW-complexes, $\widehat{g_\nu}$ induces a homeomorphism from $c_\nu$ to $\widehat{g_\nu}(c_\nu)$ equipped with the relative topology. The map $\widehat{g_{\nu}}: \globM(\mathbf{D}^{n_\nu})\to X$ is called the \textit{attaching map} of the globular cell $c_\nu$. The state $\widehat{g_\nu}(0)\in X^0$ ($\widehat{g_\nu}(1)\in X^0$ resp.)  is called the \textit{initial (final resp.) state} of $c_\nu$ and is denoted by $c_\nu^-$ ($c_\nu^+$ resp.). The integer $n_\nu+1$ is called the \textit{dimension} of the globular cell $c_\nu$. It is denoted by $\dim c_\nu$. The states of $X^0$ are also called the \textit{globular cells of dimension $0$}. By convention, a state of $X^0$ viewed as a globular cell of dimension $0$ is equal to its initial state and to its final state. Thus, for $\alpha\in X^0$, one has $\alpha=\alpha^+=\alpha^-$. The set of globular cells of the cellular decomposition $\widetilde{X}$ of $X$ is denoted by $\mathcal{C}(\widetilde{X})$. Let $x\in |X|$. There exists a unique globular cell $\dt(x)$ of the cellular decomposition $\widetilde{X}$ such that $x\in \dt(x)$ since there is the equality of sets
\[
|X| = \coprod_{c\in \mathcal{C}(\widetilde{X})} c.
\]

\begin{prop} \label{prop:glnat-branching}
	Let $X$ be a cellular multipointed $d$-space equipped with a cellular decomposition $\widetilde{X}:\lambda\to\ptop{\mathcal{M}}$. Let $\alpha\in X^0$. Every element of $\mathcal{P}^-_\alpha(X)$ is of the form $\widehat{g}\delta_z\psi$ where $\widehat{g}:\globM(\mathbf{D}^n)\to X$ is an attaching map appearing in the cellular decomposition $\widetilde{X}$ such that $\widehat{g}(0)=\alpha$ and where $z\in \mathbf{D}^n$ and $\psi\in \mathcal{I}^-$. The choice of the attaching map $\widehat{g}$ and of $\psi\in \mathcal{I}^-$ is unique if we impose the condition $z\in \mathbf{D}^n\backslash \mathbf{S}^{n-1}$. 
\end{prop}

\bpf
By \cite[Theorem~4.9]{GlobularNaturalSystem}, $X$ being cellular by hypothesis, every directed path of $\cont(X)$ starting from $\alpha$ is of the form $\gamma\phi$ where $\gamma$ is an execution path of $X$ starting from $\alpha$ and $\phi\in \mathcal{I}$. By \cite[Theorem~4.7]{GlobularNaturalSystem} which is in fact a reformulation of \cite[Theorem~6]{Moore3}, the execution path $\gamma$ can be written in a unique way as a composite of the form $\natgl(\gamma)\phi'$ with $\natgl(\gamma) = (\widehat{g_{\nu_1}}\delta_{z_1})*\dots * (\widehat{g_{\nu_n}}\delta_{z_n})$, $n\geq 1$, $\nu_i<\lambda$ and $z_i\in \mathbf{D}^{n_{\nu_i}}\backslash \mathbf{S}^{n_{\nu_i}-1}$ for $1\leq i\leq n$, and $\phi':[0,1]\to [0,n]$ a surjective non-decreasing map. When moreover $\gamma\psi\in \mathcal{P}^-_\alpha(X)$, the condition $\gamma\psi(]0,1[)\cap X^0=\varnothing$ implies that one can suppose that $n=1$. Thus $\gamma\psi=\widehat{g_{\nu_1}}\delta_{z_1}\phi'\phi$ with $\phi'\phi(]0,1[) \cap \{0,1\}=\varnothing$, which means that $\phi'\phi\in \mathcal{I}^-$. By \cite[Proposition~18]{Moore3}, the directed path $\widehat{g}\delta_z$ is regular. Thus $\psi=\phi'\phi$ is unique by \cite[Proposition~19]{Moore3}. 
\epf

\begin{lem} (Cutler's observation) \label{lem:open-delta-k}
	Let $M$ be a general topological space. If $U$ is an open subset of $M$, then the topology of $k_\Delta(U)$ is the relative topology with respect to the inclusion into $k_\Delta(M)$.
\end{lem}

\bpf
Let $V=\phi_M^{-1}(U)$. Then $V$ is an open subset of $k_\Delta(M)$, $\phi_M$ being continuous. Thus, equipped with the relative topology, $V$ is $\Delta$-generated. Hence the map $V\to U$ factors as a composite of continuous maps $V\to k_\Delta(U) \to U$. The map $V\to k_\Delta(U)$ is also bijective. This map has an inverse $k_\Delta(U)\to V$ which is obtained by factoring the map $k_\Delta(U)\to k_\Delta(M)$ induced by the inclusion $U\subset M$. Thus there is a homeomorphism $V\iso k_\Delta(U)$. 
\epf

\begin{prop} \label{prop:P-closedinclusion}
	Let $n\geq 0$. Consider a pushout diagram of cellular multipointed $d$-spaces
	\[
	\begin{tikzcd}[row sep=3em,column sep=3em]
		\globM(\mathbf{S}^{n-1}) \arrow[r,"g"] \arrow[d,hookrightarrow] & A \arrow[d,hookrightarrow,"f"] \\
		\globM(\mathbf{D}^{n}) \arrow[r,"\widehat{g}"] & \cocartesian X
	\end{tikzcd}
	\]
	Then the mapping $\gamma\mapsto f\gamma$ induces a continuous map $\mathcal{P}^-_\alpha(A) \to \mathcal{P}^-_\alpha(X)$ for all $\alpha\in A^0=X^0$ which is the identity for $\alpha\neq \widehat{g}(0)$ and a closed inclusion for $\alpha= \widehat{g}(0)$. Consequently, the mapping $\gamma\mapsto f\gamma$ induces a continuous map $\mathcal{P}^-(A) \to \mathcal{P}^-(X)$ which is a closed inclusion.
\end{prop}

\bpf
By Proposition~\ref{prop:noncontracting1}, the mapping $\gamma\mapsto f\gamma$ induces a continuous map $\mathcal{P}^-_\alpha(A) \to \mathcal{P}^-_\alpha(X)$ for all $\alpha\in A^0=X^0$. Let $\alpha\in A^0=X^0$. Let $Z\subset |X|$ and $\mathcal{P}^-_\alpha(X,Z)_{co} = \mathcal{P}^-_\alpha(X) \cap \ttop_{co}([0,1],Z)$ equipped with the relative topology for the inclusion $\mathcal{P}^-_\alpha(X,Z)_{co}\subset \ttop_{co}([0,1],Z)$. Let $\mathcal{P}^-_\alpha(X,Z)= k_\Delta(\mathcal{P}^-_\alpha(X,Z)_{co})$. Note that all these spaces are Hausdorff, $|X|$ being Hausdorff. By definition, there is the homeomorphism $\mathcal{P}^-_\alpha(X)\iso\mathcal{P}^-_\alpha(X,|X|)$. By Proposition~\ref{prop:glnat-branching}, there is the homeomorphism $\mathcal{P}^-_\alpha(A)\iso\mathcal{P}^-_\alpha(X,|A|)$. Since the functor $Z\mapsto |\globM(Z)|$ is a left Quillen adjoint for the q-model structures by \cite[Proposition~A.1]{MultipointedSubdivision}, the map $|f|:|A|\to |X|$ is a q-cofibration of spaces, and therefore a h-cofibration. Thus, $|A|$ is a closed subset of $|X|$. By \cite[Theorem~2]{vstrom1}, there exists an open subset $U$ of $|X|$ such that $|A|\subset U$ and such that $|A|$ is a retract of $U$. We now adapt Cutler's argument presented in \cite{Tyrone} to the context of the proof. It is divided in three steps.
\begin{enumerate}[leftmargin=*]
	\item We obtain the commutative diagram of spaces 
	\[
	\begin{tikzcd}[column sep=3em,row sep=0em]
		\mathcal{P}^-_\alpha(A) \arrow[r,"i"] \arrow[rr,bend right=15pt,equal]& \mathcal{P}^-_\alpha(X,U) \arrow[r,"r"] & \mathcal{P}^-_\alpha(A)
	\end{tikzcd}
	\]
	Since $i$ is the equalizer (resp. $r$ is the coequalizer) of the pair of maps $(ir,\id_{\mathcal{P}^-_\alpha(X,U)})$, there are the homeomorphisms \[\mathcal{P}^-_\alpha(A)\iso \big\{f\in \mathcal{P}^-_\alpha(X,U),f = i(r(f))\big\}\iso i(\mathcal{P}^-_\alpha(A)).\] We deduce that $\mathcal{P}^-_\alpha(A)\iso i(\mathcal{P}^-_\alpha(A))$ is a closed subset of $\mathcal{P}^-_\alpha(X,U)$, the diagonal of $\mathcal{P}^-_\alpha(X,U)$ being closed. We have proved that the map $\mathcal{P}^-_\alpha(A)\to \mathcal{P}^-_\alpha(X,U)$ is a closed inclusion of $\Delta$-generated spaces. For a similar reason, $|A|\subset U$ is also a closed inclusion of $\Delta$-generated spaces.
	\item Let $(\gamma_n)_{n\geq 0}$ be a sequence of $\mathcal{P}^-_\alpha(A)$ which converges to $\gamma_\infty$ in $\mathcal{P}^-_\alpha(X)$. Since $|A|$ is closed in $|X|$ which is Hausdorff, and since $(\gamma_n)_{n\geq 0}$ converges pointwise to $\gamma_\infty\in \ttop([0,1],|X|)$, one has $\gamma_\infty([0,1])\subset |A|$. Since $|A|\subset U$ is a closed inclusion and since $U$ is open in $|X|$, the composite map $|A|\to U \to |X|$ is a $\Delta$-inclusion, $U$ being equipped with the relative topology. This implies that $\gamma_\infty\in\ttop([0,1],|A|)$. Since $\mathcal{P}^-_\alpha(X,U)_{co}$ is an open subset of $\mathcal{P}^-_\alpha(X)_{co}$ by definition of the compact-open topology, $\mathcal{P}^-_\alpha(X,U)$ is an open subspace of $\mathcal{P}^-_\alpha(X)$ by Lemma~\ref{lem:open-delta-k}. Consequently, there exists $N\geq 0$ such that $\gamma_n\in \mathcal{P}^-_\alpha(X,U)$ for all $n\geq N$ and, moreover, the sequence $(\gamma_n)_{n\geq N}$ is convergent in $\mathcal{P}^-_\alpha(X,U)$.
	\item Putting the above facts together, we deduce that $(\gamma_n)_{n\geq 0}$ converges to $\gamma_\infty$ in $\mathcal{P}^-_\alpha(A)$. By Theorem~\ref{thm:recognizing-closed-inclusions-1}, we deduce that the map $\mathcal{P}^-_\alpha(A) \subset \mathcal{P}^-_\alpha(X)$ is a closed inclusion of $\Delta$-generated spaces. 
\end{enumerate}
Finally, for $\alpha\neq \widehat{g}(0)$, Proposition~\ref{prop:glnat-branching} implies that the closed inclusion $\mathcal{P}^-_\alpha(A) \subset \mathcal{P}^-_\alpha(X)$ is onto. Hence it is a homeomorphism in this case. 
\epf

\begin{prop} \label{prop:P-colim}
	Let $X$ be a cellular multipointed $d$-space. Let $\widehat{X}:\lambda\to\ptop{\mathcal{M}}$ be a cellular decomposition of $X$. Let $\nu\leq \lambda$ be a limit ordinal. The canonical continuous map 
	\[
	\Phi_\nu:\liminj_{\mu<\nu} \mathcal{P}^-(\widetilde{X}_\mu) \longrightarrow \mathcal{P}^-(\widetilde{X}_\nu)
	\]
	is a homeomorphism.
\end{prop}

\bpf
The attaching maps of $\widetilde{X}_\nu$ consist of the union of all attaching maps of all $\widetilde{X}_\mu$ for $\mu<\nu$ since $\widehat{X}:\lambda\to\ptop{\mathcal{M}}$ is colimit-preserving, being a cellular decomposition. Thus by Proposition~\ref{prop:glnat-branching}, the map $\Phi_\nu$ is a continuous bijection. Consider a set map $f:[0,1]\to \liminj_{\mu<\nu} \mathcal{P}^-(\widetilde{X}_\mu)$ such that the composite set map $\Phi_\nu f$ is continuous. It gives rise to a continuous map $\widehat{f}:[0,1]\p [0,1]\to \widetilde{X}_\nu$. Since the tower $\widehat{X}$ consists of q-cofibrations which are closed $T_1$-inclusions, the map $\widehat{f}:[0,1]\p [0,1]\to \widetilde{X}_\nu$ factors as a composite $\widehat{f}:[0,1]\p [0,1]\to \widetilde{X}_\mu\to\widetilde{X}_\nu$ for some $\mu<\nu$ by \cite[Proposition~2.4.2]{MR99h:55031}. By adjunction, we obtain that $\Phi_\nu f$ factors as a composite of continuous maps $[0,1]\to \mathcal{P}^-(\widetilde{X}_\nu)\cap \ttop([0,1],|X_\mu|)= \mathcal{P}^-(\widetilde{X}_\mu) \to \mathcal{P}^-(\widetilde{X}_\nu)$. Hence $f$ is continuous and the continuous bijection $\Phi_\nu$ is a $\Delta$-inclusion. Thus, it is a homeomorphism.
\epf

\begin{prop} \label{prop:preparation}
	Let $n\geq 0$. Consider a pushout diagram of cellular multipointed $d$-spaces
	\[
	\begin{tikzcd}[row sep=3em,column sep=3em]
		\globM(\mathbf{S}^{n-1}) \arrow[r,"g"] \arrow[d,hookrightarrow] & A \arrow[d,hookrightarrow,"f"] \\
		\globM(\mathbf{D}^{n}) \arrow[r,"\widehat{g}"] & \cocartesian X
	\end{tikzcd}
	\]
	Then the mapping $\gamma\mapsto f\gamma$ induces a continuous map \[\mathcal{G}^-_\alpha(A) \longrightarrow \mathcal{G}^-_\alpha(X)\] for all $\alpha\in A^0=X^0$ which is the identity for $\alpha\neq \widehat{g}(0)$ and a closed inclusion for $\alpha= \widehat{g}(0)$. Consequently, the mapping $\gamma\mapsto f\gamma$ induces a continuous map $\mathcal{G}^-(A) \to \mathcal{G}^-(X)$ which is a closed inclusion.
\end{prop}

\bpf
By Proposition~\ref{prop:noncontracting1}, the mapping $\gamma\mapsto f\gamma$ induces a continuous map $\mathcal{G}^-_\alpha(A) \to \mathcal{G}^-_\alpha(X)$ for all $\alpha\in A^0=X^0$. The case $\alpha\neq \widehat{g}(0)$ is a consequence of $\mathcal{P}^-_\alpha(A) = \mathcal{P}^-_\alpha(X)$. Assume that $\alpha=\widehat{g}(0)$. Consider the commutative diagram of $\Delta$-generated spaces 
\[
\begin{tikzcd}[column sep=3em,row sep=3em]
	\mathcal{P}_{\widehat{g}(0)}^-(A) \arrow[r,"p"] \arrow[d,"f"']& \mathcal{G}^-_{\widehat{g}(0)}(A) \arrow[d,"\overline{f}"] \\
	\mathcal{P}_{\widehat{g}(0)}^-(X) \arrow[r,"q"] & \mathcal{G}^-_{\widehat{g}(0)}(X)
\end{tikzcd}
\]
By Proposition~\ref{prop:P-closedinclusion}, the map $f$ is a closed inclusion. Let $\gamma\in q^{-1}(\overline{f}(\mathcal{G}^-_{\widehat{g}(0)}(A)))$. Then $\gamma=\widehat{k}\delta_z\phi$ where $\widehat{k}:\globM(\mathbf{D}^N)\to X$ is an attaching map of the cellular decomposition of $A$, $z\in \mathbf{D}^N\backslash\mathbf{S}^{N-1}$ and $\phi\in \mathcal{I}^-$. We deduce that $\gamma\in f(\mathcal{P}_{\widehat{g}(0)}^-(A))$. Thanks to the inclusion $q^{-1}(\overline{f}(\mathcal{G}^-_{\widehat{g}(0)}(A)))\subset f(\mathcal{P}_{\widehat{g}(0)}^-(A))$, the proof is complete by Corollary~\ref{cor:recognizing-closed-inclusions-3}.
\epf

\begin{prop} \label{prop:G-colim}
	Let $X$ be a cellular multipointed $d$-space. Let $\widehat{X}:\lambda\to\ptop{\mathcal{M}}$ be a cellular decomposition of $X$. Let $\nu\leq \lambda$ be a limit ordinal. The canonical continuous map 
	\[
	\Psi_\nu:\liminj_{\mu<\nu} \mathcal{G}^-(\widetilde{X}_\mu) \longrightarrow \mathcal{G}^-(\widetilde{X}_\nu)
	\]
	is a homeomorphism.
\end{prop}

\bpf
Consider the commutative diagram of spaces 
\[
\begin{tikzcd}[row sep=3em,column sep=3em]
	\liminj\limits_{\mu<\nu} \mathcal{P}^-(\widetilde{X}_\mu) \arrow[r,"p"] \arrow[d,"\iso"']& \liminj\limits_{\mu<\nu} \mathcal{G}^-(\widetilde{X}_\mu) \arrow[d] \\
	 \mathcal{P}^-(\widetilde{X}_\nu) \arrow[r] & \mathcal{G}^-(\widetilde{X}_\nu)=\mathcal{P}^-(\widetilde{X}_\nu)/\!\simeq^-
\end{tikzcd}
\]
The left vertical map is a homeomorphism by Proposition~\ref{prop:P-colim}. The map $p$ is onto. The right vertical map is a continuous bijection by Proposition~\ref{prop:glnat-branching}. We conclude that the right vertical map is a homeomorphism using Corollary~\ref{cor:recognizing-closed-inclusions-4}.
\epf

\begin{lem} \label{lem:converge} \cite[Lemma~2]{Moore3}
	Let $X$ be a sequential topological space. Let $x_\infty\in X$. Let $(x_n)_{n\geq 0}$ be a sequence such that $x_\infty$ is a limit point of all subsequences. Then the sequence $(x_n)_{n\geq 0}$ converges to $x_\infty$.
\end{lem}

\bth \label{thm:germ-glob}
Let $n\geq 0$. Consider a pushout diagram of cellular multipointed $d$-spaces
\[
\begin{tikzcd}[row sep=3em,column sep=3em]
	\globM(\mathbf{S}^{n-1}) \arrow[r,"g"] \arrow[d,hookrightarrow] & A \arrow[d,hookrightarrow,"f"] \\
	\globM(\mathbf{D}^{n}) \arrow[r,"\widehat{g}"] & \cocartesian X
\end{tikzcd}
\]
Then the map $f:A\to X$ induces a bijection $A^0\iso X^0$ and for all $\alpha\in A^0$ a continuous map $\mathcal{G}^-_\alpha(A) \to \mathcal{G}^-_\alpha(X)$. Moreover there is the pushout diagram of spaces 
\[
\begin{tikzcd}[row sep=3em,column sep=3em]
	\mathbf{S}^{n-1} \arrow[r] \arrow[d,hookrightarrow] & \mathcal{G}^-_{g(0)}(A) \arrow[d,hookrightarrow] \\
	\mathbf{D}^{n} \arrow[r] & \cocartesian \mathcal{G}^-_{g(0)}(X)
\end{tikzcd}
\]
and for all $\alpha\in A^0\backslash \{g(0)\}$ the homeomorphism $\mathcal{G}^-_{\alpha}(A)\iso \mathcal{G}^-_{\alpha}(X)$.
\eth

\bpf The case $\alpha\in A^0\backslash \{g(0)\}$ is treated in Proposition~\ref{prop:preparation}. Assume that $\alpha=g(0)=\widehat{g}(0)$. The multipointed $d$-space $A$ being cellular, consider a cellular decomposition $\widetilde{A}:\lambda\to \ptop{\mathcal{M}}$ of $A$ such that \begin{itemize}
	\item The multipointed $d$-space $\widetilde{A}_0$ is $(A^0,A^0,\varnothing)$.
	\item For all $\nu<\lambda$, there is a pushout diagram of multipointed $d$-spaces 
	\[\begin{tikzcd}[row sep=3em, column sep=3em]
		\globM(\mathbf{S}^{n_\nu-1}) \arrow[d,rightarrowtail] \arrow[r,"g_\nu"] & \widetilde{A}_\nu \arrow[d,rightarrowtail] \\
		\globM(\mathbf{D}^{n_\nu}) \arrow[r,"\widehat{g_\nu}"] & \cocartesian \widetilde{A}_{\nu+1}
	\end{tikzcd}\]
	with $n_\nu \geq 0$. 
\end{itemize}
By adding the cell corresponding to the pushout diagram of the statement of the theorem, we obtain a cellular decomposition $\widetilde{X}:\lambda+1\to \ptop{\mathcal{M}}$ of $X$ with $\widetilde{A}_\nu=\widetilde{X}_\nu$ for $\nu\leq \lambda$. Since the underlying space functor is a left adjoint by \cite[Proposition~A.2]{MultipointedSubdivision}, it is colimit-preserving. Thus, there is the pushout diagram of spaces 
\[
\begin{tikzcd}[row sep=3em,column sep=3em]
	{|\globM(\mathbf{S}^{n-1})|} \arrow[r,"|g|"] \arrow[d] & {|A|} \arrow[d,"|f|"] \\
	{|\globM(\mathbf{D}^{n})|} \arrow[r,"|\widehat{g}|"] & \cocartesian |X|
\end{tikzcd}
\]
Consequently the continuous map $|A|\to |X|$ is one-to-one. Using Corollary~\ref{cor:globe-branching} and Proposition~\ref{prop:noncontracting1}, we obtain the commutative diagram of topological spaces 
\[
\begin{tikzcd}[row sep=3em,column sep=3em]
	\mathbf{S}^{n-1} \arrow[r] \arrow[d,hookrightarrow] & \mathcal{G}^-_{g(0)}(A) \arrow[d,hookrightarrow] \arrow[ddr,bend left=20pt,"k_1"]\\
	\mathbf{D}^{n} \arrow[r] \arrow[rrd,bend right=20pt,"k_2"']&  \cocartesian Z \arrow[rd,dashed,"\exists ! k"]\\ 
	&& \mathcal{G}^-_{g(0)}(X)
\end{tikzcd}
\]
By Proposition~\ref{prop:glnat-branching}, every element of $\mathcal{G}^-_{g(0)}(X)$ is the equivalence class of a directed path either of the form $\widehat{g}\delta_z\phi$ for $z\in \mathbf{D}^{n}\backslash \mathbf{S}^{n-1}$ and $\phi\in \mathcal{I}^-$ or of the form $\widehat{g'}\delta_{z'}\phi'$ where $\widehat{g'}$ is an attaching map of the cellular decomposition of $A$ and $\phi'\in\mathcal{I}^-$. Thus the continuous map $k:Z \to \mathcal{G}^-_{g(0)}(X)$ is surjective. The uniqueness of Proposition~\ref{prop:glnat-branching} implies that the continuous map $k:Z \to \mathcal{G}^-_{g(0)}(X)$ is one-to-one. Consequently, $k$ is a continuous bijection. Consider a sequence $(\tr{\widehat{g_k}\delta_{z_k}\phi_k}^-)_{k\geq 0}$ of $\mathcal{G}^-_{g(0)}(X)$ converging to $L$. There are two mutually exclusive cases:
\begin{enumerate}
	\item $\{k\geq 0 \mid \widehat{g_k}=\widehat{g}\}$ infinite: by extracting a subsequence, we can suppose that the sequence is of the form $(\tr{\widehat{g}\delta_{z_k}}^-)_{k\geq 0}$ with $z_k\in \mathbf{D}^n$; in this case, one has $\tr{\widehat{g}\delta_{z_k}}^-=k_2(z_k)$. By extracting a subsequence, we can suppose that the sequence $(z_k)_{k\geq 0}$ of the compact metrizable space $\mathbf{D}^n$ converges to $z_\infty$; thus the sequence $(\tr{\widehat{g_k}\delta_{z_k}\phi_k}^-)_{k\geq 0}$ of $Z$ has a limit point which is necessarily $L$.
	\item $\{k\geq 0 \mid \widehat{g_k}=\widehat{g}\}$ finite: by extracting a subsequence, we can suppose that the sequence $(\tr{\widehat{g_k}\delta_{z_k}\phi_k}^-)_{k\geq 0}$ belongs to $k_1(\mathcal{G}^-_{g(0)}(A))$; the map $k_1$ being a closed inclusion by Proposition~\ref{prop:preparation}, the sequence $(\tr{\widehat{g_k}\delta_{z_k}\phi_k}^-)_{k\geq 0}$ converges to $L$ in $\mathcal{G}^-_{g(0)}(A)$; thus it converges in $Z$ as well; we have proved that the sequence $(\tr{\widehat{g_k}\delta_{z_k}\phi_k}^-)_{k\geq 0}$ of $Z$ has a limit point which is necessarily $L$. 
\end{enumerate}
We have proved that every subsequence of $(\tr{\widehat{g_k}\delta_{z_k}\phi_k}^-)_{k\geq 0}$ in $Z$ has a limit point which is necessarily $L$. By Lemma~\ref{lem:converge}, $Z$ being sequential, the sequence $(\tr{\widehat{g_k}\delta_{z_k}\phi_k}^-)_{k\geq 0}$ converges to $L$ in $Z$. We deduce that $Z$ and $\mathcal{G}^-_{g(0)}(X)$  have the same convergent sequences. Since both spaces are sequential, being $\Delta$-generated, we conclude that the continuous bijection $k:Z\to \mathcal{G}^-_{g(0)}(X)$ is a homeomorphism. 
\epf

We cannot use Theorem~\ref{thm:recognizing-closed-inclusions-1} to prove that $k:Z\to \mathcal{G}^-_{g(0)}(X)$ is a closed inclusion (and hence a homeomorphism) since we do not know yet that $\mathcal{G}^-_{g(0)}(X)$, which is equipped with a final topology, is $\Delta$-Hausdorff.

\begin{cor} \label{cor:push-glob}
Let $n\geq 0$. Consider a pushout diagram of cellular multipointed $d$-spaces
\[
\begin{tikzcd}[row sep=3em,column sep=3em]
	\globM(\mathbf{S}^{n-1}) \arrow[r,"g"] \arrow[d,hookrightarrow] & A \arrow[d,hookrightarrow,"f"] \\
	\globM(\mathbf{D}^{n}) \arrow[r,"\widehat{g}"] & \cocartesian X
\end{tikzcd}
\]
Then the map $f:A\to X$ induces a continuous map $\mathcal{G}^-(A) \to \mathcal{G}^-(X)$. Moreover there is the pushout diagram of spaces 
\[
\begin{tikzcd}[row sep=3em,column sep=3em]
	\mathbf{S}^{n-1} \arrow[r] \arrow[d,hookrightarrow] & \mathcal{G}^-(A) \arrow[d,hookrightarrow] \\
	\mathbf{D}^{n} \arrow[r] & \cocartesian \mathcal{G}^-(X)
\end{tikzcd}
\]
\end{cor}

\begin{prop} \label{prop:push-closed-incl} \cite[Proposition~B.15]{leftproperflow} Consider the pushout diagram in the category of $\Delta$-generated spaces 
\[
\begin{tikzcd}[row sep=3em,column sep=3em]
	A \arrow[r,"h"] \arrow[d,"f"'] & X \arrow[d,"g"] \\
	B \arrow[r,"k"]& \cocartesian Y
\end{tikzcd}
\]
If the map $f:A\to B$ is a closed inclusion, then the map $g:X\to Y$ is a closed inclusion. If moreover $B$ and $X$ are $\Delta$-Hausdorff, then $Y$ is $\Delta$-Hausdorff. 
\end{prop}

\begin{prop} \label{prop:alreadyH} \cite[Proposition~B.16]{leftproperflow} Consider a transfinite tower $X$ of $\Delta$-Hausdorff $\Delta$-generated spaces such that for every ordinal $\alpha$, the map $X_\alpha\to X_{\alpha+1}$ is one-to-one. Then the colimit $\liminj X$ equipped with the final topology is $\Delta$-Hausdorff.
\end{prop}

\begin{cor} \label{cor:germ-glob}
Let $X$ be a q-cofibrant multipointed $d$-space. The final quotient $\mathcal{G}^-_\alpha(X)$ is q-cofibrant. In particular, it is $\Delta$-Hausdorff, which implies that the equivalence relation $\simeq^-$ has a closed graph in this case. 
\end{cor}

\bpf
Any q-cofibrant multipointed $d$-space being a retract of a cellular multipointed $d$-space, and the retract of a $\Delta$-Hausdorff space being $\Delta$-Hausdorff, one can assume that $X$ is cellular. We consider a cellular decomposition $\widetilde{X}:\lambda\to \ptop{\mathcal{M}}$ of $X$ such that \begin{itemize}
	\item The multipointed $d$-space $\widetilde{X}_0$ is $(X^0,X^0,\varnothing)$.
	\item For all $\nu<\lambda$, there is a pushout diagram of multipointed $d$-spaces 
	\[\begin{tikzcd}[row sep=3em, column sep=3em]
		\globM(\mathbf{S}^{n_\nu-1}) \arrow[d,rightarrowtail] \arrow[r,"g_\nu"] & \widetilde{X}_\nu \arrow[d,rightarrowtail] \\
		\globM(\mathbf{D}^{n_\nu}) \arrow[r,"\widehat{g_\nu}"] & \cocartesian \widetilde{X}_{\nu+1}
	\end{tikzcd}\]
	with $n_\nu \geq 0$. 
\end{itemize}
We prove by a transfinite induction on $\nu\leq \lambda$ that $\mathcal{G}^-_\alpha(\widetilde{X}_\nu)$ is q-cofibrant for all $\alpha\in X^0$. For all $\alpha\in X^0$, $\mathcal{G}^-_\alpha(\widetilde{X}_0)=\varnothing$, which is $\Delta$-Hausdorff. The passage from $\nu<\lambda$ to $\nu+1$ is a consequence of Theorem~\ref{thm:germ-glob} and Proposition~\ref{prop:push-closed-incl}. The case $\nu$ limit ordinal is a consequence of Proposition~\ref{prop:G-colim} and Proposition~\ref{prop:alreadyH}.
\epf

\section{Comparing multipointed d-spaces and flows}
\label{section:comparing-mdtop-flow}

We need to make some additional reminders about the functor $\dcat:\ptop{\mathcal{M}}\to \dtop$ before proving Theorem~\ref{thm:dgerm-cgerm}.

\begin{prop}\label{prop:colim} \cite[Proposition~7.1]{GlobularNaturalSystem} Consider a pushout diagram of multipointed $d$-spaces of the form 
	\[
	\begin{tikzcd}[row sep=3em, column sep=3em]
		\globM(\mathbf{S}^{n-1}) \arrow[r] \arrow[d] & A \arrow[d] \\
		\globM(\mathbf{D}^{n}) \arrow[r] & \cocartesian B
	\end{tikzcd}
	\]
	with $A$ cellular and $n\geq 0$. Then there is a pushout diagram of flows 
	\[
	\begin{tikzcd}[row sep=3em, column sep=3em]
		\glob(\mathbf{S}^{n-1}) \arrow[r] \arrow[d] & \dcat(A) \arrow[d] \\
		\glob(\mathbf{D}^{n}) \arrow[r] & \cocartesian \dcat(B).
	\end{tikzcd}
	\]
\end{prop}

\begin{prop} \label{prop:almost-accessible} \cite[Corollary~4.3]{MultipointedSubdivision}
	Let $X:\lambda\to\ptop{\mathcal{M}}$ be a transfinite tower of q-cofibrations between q-cofibrant multipointed $d$-spaces. Then the canonical map \[\liminj_{\nu<\lambda} \dcat(X_\nu)\longrightarrow \dcat(\liminj_{\nu<\lambda} X_\nu)\] is an isomorphism of flows.
\end{prop}

Theorem~\ref{thm:dgerm-cgerm} relates a definition of the branching space in the semicategorical framework of flows to a geometric interpretation in terms of germs of short directed paths. This geometric interpretation is used in Proposition~\ref{prop:short-branching}, Proposition~\ref{prop:noncontracting2}, and then in Theorem~\ref{thm:noncontracting3}.

\bth \label{thm:dgerm-cgerm}
Let $X$ be a q-cofibrant multipointed $d$-space. For all $\alpha\in X^0$, there is a homeomorphism 
\[
\mathcal{G}^-_\alpha(X) \iso \P^-_\alpha\dcat(X) \simeq \hop^-_\alpha \dcat(X).
\]
\eth

\bpf
By \cite[Theorem~15]{Moore3}, the flow $\dcat(X)$ is q-cofibrant, $X$ being q-cofibrant by hypothesis. Thus for any q-cofibrant replacement functor $(-)^{cof}$ of the q-model category of flows, the natural map $\dcat(X)^{cof}\to \dcat(X)$ is a trivial fibration between q-cofibrant flows. Hence, the natural map $\hop^-\dcat(X)\to \P^-\dcat(X)$ is a weak homotopy equivalence between q-cofibrant spaces, the functor $\P^-:\dtop\to \top$ being a left Quillen functor, and therefore it is a homotopy equivalence. Any q-cofibrant multipointed $d$-space being a retract of a cellular multipointed $d$-space, and the retract of a homeomorphism being a homeomorphism, one can assume that $X$ is cellular to prove the homeomorphism $\mathcal{G}^-_\alpha(X) \iso \P^-_\alpha\dcat(X)$. We consider a cellular decomposition $\widetilde{X}:\lambda\to \ptop{\mathcal{M}}$ of $X$ such that \begin{itemize}
	\item The multipointed $d$-space $\widetilde{X}_0$ is $(X^0,X^0,\varnothing)$.
	\item For all $\nu<\lambda$, there is a pushout diagram of multipointed $d$-spaces 
	\[\begin{tikzcd}[row sep=3em, column sep=3em]
		\globM(\mathbf{S}^{n_\nu-1}) \arrow[d,rightarrowtail] \arrow[r,"g_\nu"] & \widetilde{X}_\nu \arrow[d,rightarrowtail] \\
		\globM(\mathbf{D}^{n_\nu}) \arrow[r,"\widehat{g_\nu}"] & \cocartesian \widetilde{X}_{\nu+1}
	\end{tikzcd}\]
	with $n_\nu \geq 0$. 
\end{itemize}
We prove by a transfinite induction on $\nu\leq \lambda$ the homeomorphism $\mathcal{G}^-_\alpha(\widetilde{X}_\nu) \iso \P^-_\alpha\dcat(\widetilde{X}_\nu)$ as follows. The case $\nu=0$ comes from the fact that $\mathcal{G}^-_\alpha(\widetilde{X}_0) \iso \P^-_\alpha\dcat(\widetilde{X}_0)=\varnothing$. Assume that there is a homeomorphism $\mathcal{G}^-_\alpha(\widetilde{X}_\nu) \iso \P^-_\alpha\dcat(\widetilde{X}_\nu)$ for some $\nu<\lambda$. Using Corollary~\ref{cor:push-glob}, we obtain the following pushout diagram of spaces 
\[\begin{tikzcd}[row sep=3em, column sep=3em]
	\mathbf{S}^{n_\nu-1} \arrow[d,rightarrowtail] \arrow[r] & \mathcal{G}^-(\widetilde{X}_\nu) \arrow[d,rightarrowtail] \\
	\mathbf{D}^{n_\nu} \arrow[r] & \cocartesian \mathcal{G}^-(\widetilde{X}_{\nu+1})
\end{tikzcd}\]
Using Proposition~\ref{prop:colim} and the fact that the functor $\P^-:\dtop\to \top$ is a left adjoint by Theorem~\ref{thm:def-Cminus}, we obtain the 
pushout diagram of spaces 
\[\begin{tikzcd}[row sep=3em, column sep=3em]
	\mathbf{S}^{n_\nu-1} \arrow[d,rightarrowtail] \arrow[r] & \P^-(\dcat(\widetilde{X}_\nu)) \arrow[d,rightarrowtail] \\
	\mathbf{D}^{n_\nu} \arrow[r] & \cocartesian \P^-(\dcat(\widetilde{X}_{\nu+1}))
\end{tikzcd}\]
Hence the passage from $\nu$ to $\nu+1$. Assume now that $\nu\leq\lambda$ is a limit ordinal. One obtains 
\begin{align*}
	\mathcal{G}^-(\widetilde{X}_\nu) &\iso \mathcal{G}^-(\liminj_{\mu<\nu}\widetilde{X}_\mu) \\
	& \iso \liminj_{\mu<\nu} \mathcal{G}^-(\widetilde{X}_\mu)\\
	&\iso \liminj_{\mu<\nu} \P^-\dcat(\widetilde{X}_\mu) \\
	&\iso \P^-\bigg(\liminj_{\mu<\nu} \dcat(\widetilde{X}_\mu)\bigg) \\
	&\iso \P^-\bigg( \dcat(\liminj_{\mu<\nu}\widetilde{X}_\mu)\bigg) \\
	&\iso \P^- \dcat(\widetilde{X}_\nu),
\end{align*}
the first and sixth homeomorphisms since the tower $\widetilde{X}:\lambda\to \ptop{\mathcal{M}}$ is colimit-preserving, being a cellular decomposition, the second homeomorphism by Proposition~\ref{prop:G-colim}, the third homeomorphism by induction hypothesis, the fourth homeomorphism since $\P^-$ is a left adjoint by Theorem~\ref{thm:def-Cminus} and finally the fifth homeomorphism by Proposition~\ref{prop:almost-accessible}.
\epf

This leads us to a purely topological definition of the branching homology of a q-cofibrant multipointed $d$-space, without any use of the functor $\dcat:\ptop{\mathcal{M}}\to \dtop$.

\begin{cor} \label{cor:gooddef-branching-homology-mdtop}
	For any q-cofibrant multipointed $d$-space $X$, the branching homology $H_*^-(X)=H_*^-(\mathbf{L}\dcat(X))$ can be defined as follows: 
	\begin{enumerate}
		\item for $n\geq 1$, $H_{n+1}^-(X):=H_n(\mathcal{G}^-(X))$
		\item  $H_1^-(X):=\ker(\epsilon)/\im(\partial)$
		\item $H_0^-(X):=\mathbb{Z}[X^0]/\im(\epsilon)$
	\end{enumerate}
	with the augmentation $\epsilon:\mathbb{Z}[\mathcal{G}^-(X)]\to \mathbb{Z}[X^0]$ defined by $\epsilon(\gamma)=\gamma(0)$ and the map $\partial:\mathbb{Z}[\top([0,1],\mathcal{G}^-(X))]\to \mathbb{Z}[\mathcal{G}^-(X)]$ defined by $\partial(f)=f(0)-f(1)$.
\end{cor}

\section{Globular subdivision and branching homology}
\label{section:glb-sbd}

\bd \cite[Definition~4.10]{diCW} \cite[Definition~9.1]{GlobularNaturalSystem} \label{def:globular-sbd}
A map of multipointed $d$-spaces $f:X\to Y$ is a \textit{globular subdivision} if both $X$ and $Y$ are cellular and if $f$ induces a homeomorphism between the underlying topological spaces of $X$ and $Y$. We say that $Y$ is a \textit{globular subdivision} of $X$ when there exists such a map. This situation is denoted by \[\begin{tikzcd}[cramped]f:X\arrow[r,"\sbd"]&Y.\end{tikzcd}\]
\ed

\begin{prop} \label{prop:subd-gl-globe} \cite[Proposition~6.3]{MultipointedSubdivision}
	Let $n\geq 0$. Consider a finite set \[F\subset |\globM(\mathbf{D}^n)| \backslash |\globM(\mathbf{S}^{n-1})|.\] Then the following data assemble into a multipointed $d$-space denoted by $\globM(\mathbf{D}^n)_F$:
	\begin{itemize}
		\item The set of states is $\{0,1\} \cup F$.
		\item The underlying space is $|\globM(\mathbf{D}^n)|$.
		\item For all $\alpha\neq \beta \in \{0,1\} \cup F$, $\P^{top}_{\alpha,\beta} \globM(\mathbf{D}^n)_F = \vec{P}(\cont(\globM(\mathbf{D}^n)))(\alpha,\beta)$. 
		\item For all $\alpha\in \{0,1\} \cup F$, $\P^{top}_{\alpha,\alpha} \globM(\mathbf{D}^n)_F=\varnothing$.  
	\end{itemize}
	In particular, there is the isomorphism of multipointed $d$-spaces \[\globM(\mathbf{D}^n)\iso \globM(\mathbf{D}^n)_\varnothing.\] 
\end{prop}

We need to recall a homotopical lemma from \cite{MultipointedSubdivision} to prove Theorem~\ref{thm:elm-sbd} which could have been put in \cite{MultipointedSubdivision}.

\begin{prop} \label{prop:hocolim-tower} \cite[Proposition~2.4]{MultipointedSubdivision}
	Let $\K$ be a cocomplete category. Let $I$ be a set of maps of $\K$. Let $\lambda$ be a transfinite ordinal. Consider two transfinite towers $A:\lambda \to \C$ and $B:\lambda \to \C$ and a natural transformation $f:A\Rightarrow B$. Assume that for each ordinal $\nu<\lambda$, the map $B_\nu \sqcup_{A_\nu} A_{\nu+1} \to B_{\nu+1}$ belongs to $\cell(I)$. Then the map 
	\[
	B_0 \sqcup_{A_0} A_\lambda \longrightarrow B_\lambda
	\]
	belongs to $\cell(I)$ as well.
\end{prop}

We also need this well-known categorical fact:

\begin{prop} \label{prop:about-pushout-squares} (\cite[Proposition~2.1 and Corollary~2.2]{MultipointedSubdivision})
	Let $\K$ be a cocomplete category. Consider the commutative diagram of objects of $\K$
	\[\begin{tikzcd}[row sep=3em, column sep=3em]
		A \arrow[r]\arrow[d] \arrow[phantom,"\underline{C}",pos=0.5,rd]  & B \arrow[r] \arrow[d] \arrow[phantom,"\underline{D}",pos=0.5,rd]& C \arrow[d] \\
		D \arrow[r]& E \arrow[r] & F
	\end{tikzcd}\]
	If $\underline{C}$ and $\underline{D}$ are pushout squares, then the composite square $\underline{C}+\underline{D}$ is a pushout square. If $\underline{C}$ and $\underline{C}+\underline{D}$ are pushout squares, then the commutative square $\underline{D}$ is a pushout square.
\end{prop}

\bth \label{thm:elm-sbd}
Every globular subdivision is a transfinite composition of pushouts of globular subdivisions of the form $\globM(\mathbf{D}^{n})\to \globM(\mathbf{D}^{n})_F$.
\eth

\bpf
Consider a globular subdivision \begin{tikzcd}[cramped]f:X\arrow[r,"\sbd"]&Y\end{tikzcd}. Let $\widetilde{X}:\lambda\to \ptop{\mathcal{M}}$ be a cellular decomposition of $X$. By \cite[Theorem~7.10]{MultipointedSubdivision}, there exists a transfinite tower of cellular multipointed $d$-spaces $\widetilde{Y}:\lambda\to \ptop{\mathcal{M}}$ and a map of transfinite towers $\widetilde{X}\to \widetilde{Y}$ such that the colimit is the globular subdivision $X\to Y$ and such that for all $\nu<\lambda$, there is a commutative diagram of multipointed $d$-spaces of the form
\[\begin{tikzcd}[row sep=3em, column sep=3em]
	\globM(\mathbf{S}^{n_\nu-1}) \arrow[d,rightarrowtail] \arrow[r,"g_\nu"] \arrow[dd,bend right=60pt,rightarrowtail]& \widetilde{X}_\nu \arrow[d,rightarrowtail] \arrow[r,"\sbd"] & \widetilde{Y}_\nu \arrow[dd,rightarrowtail]\\
	\globM(\mathbf{D}^{n_\nu}) \arrow[d,"\sbd"]\arrow[r,"\widehat{g_\nu}"] & \cocartesian \widetilde{X}_{\nu+1}\arrow[rd,shorten >=0.5em,"\sbd"]  \\
	\globM(\mathbf{D}^{n_\nu})_{F(c_\nu)} \arrow[rr] & & \widetilde{Y}_{\nu+1} \arrow[lu, phantom, "\ulcorner"{font=\Large}, pos=-0.1]
\end{tikzcd}\]
where the commutative square consisting of the vertices $\globM(\mathbf{S}^{n_\nu-1})$, $\globM(\mathbf{D}^{n_\nu})_{F(c_\nu)}$, $\widetilde{Y}_\nu$ and $\widetilde{Y}_{\nu+1}$ is a pushout square. Let us rewrite the above diagram as follows: 
\[\begin{tikzcd}[row sep=3em, column sep=3em]
	\globM(\mathbf{S}^{n_\nu-1}) \arrow[d,rightarrowtail] \arrow[r,"g_\nu"] \arrow[dd,bend right=60pt,rightarrowtail]& \widetilde{X}_\nu \arrow[d,rightarrowtail] \arrow[r,"\sbd"] & \widetilde{Y}_\nu \arrow[d,rightarrowtail]\\
	\globM(\mathbf{D}^{n_\nu}) \arrow[d,"\sbd"]\arrow[r,"\widehat{g_\nu}"] & \cocartesian \widetilde{X}_{\nu+1}\arrow[r,"\sbd"] & \cocartesian \widetilde{Y}_\nu \sqcup_{\widetilde{X}_\nu} \widetilde{X}_{\nu+1} \arrow[d,"\sbd"]\\
	\globM(\mathbf{D}^{n_\nu})_{F(c_\nu)} \arrow[rr] & & \widetilde{Y}_{\nu+1}
\end{tikzcd}\]
Using Proposition~\ref{prop:about-pushout-squares}, we deduce that the commutative square consisting of the vertices $\globM(\mathbf{S}^{n_\nu-1})$, $\globM(\mathbf{D}^{n_\nu})$, $\widetilde{Y}_\nu$ and $\widetilde{Y}_\nu \sqcup_{\widetilde{X}_\nu} \widetilde{X}_{\nu+1}$ is a pushout square. Using Proposition~\ref{prop:about-pushout-squares}, we deduce that the bottom commutative square is a pushout square as well:
\[\begin{tikzcd}[row sep=3em, column sep=3em]
	\globM(\mathbf{S}^{n_\nu-1}) \arrow[d,rightarrowtail] \arrow[r,"g_\nu"] \arrow[dd,bend right=60pt,rightarrowtail]& \widetilde{X}_\nu \arrow[d,rightarrowtail] \arrow[r,"\sbd"] & \widetilde{Y}_\nu \arrow[d,rightarrowtail]\\
	\globM(\mathbf{D}^{n_\nu}) \arrow[d,"\sbd"]\arrow[r,"\widehat{g_\nu}"] & \cocartesian \widetilde{X}_{\nu+1}\arrow[r,"\sbd"] & \cocartesian \widetilde{Y}_\nu \sqcup_{\widetilde{X}_\nu} \widetilde{X}_{\nu+1} \arrow[d,"\sbd"]\\
	\globM(\mathbf{D}^{n_\nu})_{F(c_\nu)} \arrow[rr] & & \cocartesian \widetilde{Y}_{\nu+1}
\end{tikzcd}\]
We have proved that for all $\nu<\lambda$, the map of multipointed $d$-spaces 
\[
\widetilde{Y}_\nu \sqcup_{\widetilde{X}_\nu} \widetilde{X}_{\nu+1} \longrightarrow \widetilde{Y}_{\nu+1}
\]
is a pushout of a globular subdivision of the form $\globM(\mathbf{D}^{n})\to \globM(\mathbf{D}^{n})_F$. The proof is complete thanks to Proposition~\ref{prop:hocolim-tower}.
\epf

\begin{prop} \label{prop:short-branching}
	Let $X$ be a multipointed $d$-space. Let $\alpha\in X^0$. Let $\phi\in \mathcal{I}^-$. Let $\mathcal{G}^-_\alpha(X,\phi) = \mathcal{P}^-_\alpha(X,\phi)/\!\simeq^-$ equipped with the final topology with
	\[
	\mathcal{P}^-_\alpha(X,\phi) = \big\{\gamma\phi\mid \gamma\in \mathcal{P}^-_\alpha(X)\}
	\]
	Then there is a homeomorphism $\mathcal{G}^-_\alpha(X,\phi)\iso \mathcal{G}^-_\alpha(X)$.
\end{prop}

\bpf For all $\gamma\in \mathcal{P}^-_\alpha(X)$, $\gamma\phi\in \mathcal{P}^-_\alpha(X)$. We obtain a continuous inclusion $\mathcal{P}^-_\alpha(X,\phi)\subset \mathcal{P}^-_\alpha(X)$ which induces a continuous map $\Psi_1:\mathcal{G}^-_\alpha(X,\phi)\to \mathcal{G}^-_\alpha(X)$. Conversely, let $\gamma \in \mathcal{P}^-_\alpha(X)$. Then the precomposition by $\phi$ gives rise to a map $\gamma\phi\in \mathcal{P}^-_\alpha(X,\phi)$. We obtain a continuous map $\mathcal{P}^-_\alpha(X) \to \mathcal{P}^-_\alpha(X,\phi)$ which induces a continuous map $\Psi_2:\mathcal{G}^-_\alpha(X)\to \mathcal{G}^-_\alpha(X,\phi)$. One has $\Psi_2\Psi_1(\gamma\phi)=\gamma\phi\phi\simeq^-\gamma\phi$ and $\Psi_1\Psi_2(\gamma)=\gamma\phi\simeq^-\gamma$. Hence the proof is complete.
\epf

\begin{prop} \label{prop:noncontracting2}
	Let $n\geq 0$. Let $F$ be a finite subset of $|\globM(\mathbf{D}^n)|\backslash |\globM(\mathbf{S}^{n-1})|$. Consider a pushout diagram of cellular multipointed $d$-spaces
	\[
	\begin{tikzcd}[row sep=3em,column sep=3em]
		\globM(\mathbf{D}^{n}) \arrow[r,"g"] \arrow[d,"\sbd"'] & A \arrow[d,"f"] \\
		\globM(\mathbf{D}^{n})_F \arrow[r,"\widehat{g}"] & \cocartesian X
	\end{tikzcd}
	\]
	Then the mapping $\gamma\mapsto f\gamma$ induces a homeomorphism $\mathcal{G}^-_\alpha(A) \to \mathcal{G}^-_{f(\alpha)}(X)$ for all $\alpha\in A^0$. 
\end{prop}

\bpf
There is nothing to prove if $F$ is empty. Assume that $F$ is nonempty. Let $F=\{(z_1,t_1),\dots,(z_p,t_p)\}$. Since $F\subset |\globM(\mathbf{D}^n)|\backslash |\globM(\mathbf{S}^{n-1})|$, one has $F\cap \{0,1\}=\varnothing$. Thus $t_1,\dots,t_p \in ]0,1[$. Let $h=\min(t_1,\dots,t_p)/2 \in ]0,1[$. Let $\phi\in \mathcal{I}^-$ defined by $\phi(t)=ht$. Let $\gamma\in \mathcal{P}^-_\alpha(A)$. Then $f\gamma\phi(]0,1[)\cap X^0=\varnothing$ whether $\alpha=g(0)$ or not, $|f|$ being a homeomorphism. We obtain a well-defined map $\mathcal{P}^-_\alpha(A)\to \mathcal{P}^-_{f(\alpha)}(X,\phi)$ which gives rise to a continuous map $\mathcal{G}^-_\alpha(A)\to \mathcal{G}^-_{f(\alpha)}(X,\phi)\iso \mathcal{G}^-_{f(\alpha)}(X)$ by Proposition~\ref{prop:short-branching}. The underlying continuous map $|f|$ being a homeomorphism, there is an inclusion $\mathcal{P}^-_{f(\alpha)}(X,\phi)\subset \mathcal{P}^-_\alpha(A,\phi)$ which gives rise to a continuous map $\mathcal{G}^-_{f(\alpha)}(X,\phi)\iso \mathcal{G}^-_{f(\alpha)}(X)\to \mathcal{G}^-_\alpha(A,\phi)\iso \mathcal{G}^-_\alpha(A)$ by Proposition~\ref{prop:short-branching} again. Hence the proof is complete.
\epf

\bth \label{thm:noncontracting3}
Any globular subdivision \begin{tikzcd}[cramped]f:X\arrow[r,"\sbd"]&Y\end{tikzcd} induces a homeomorphism \[\mathcal{G}^-_\alpha(X) \iso \mathcal{G}^-_{f(\alpha)}(Y)\] for all $\alpha\in X^0$.
\eth

\bpf
By Theorem~\ref{thm:elm-sbd} and Proposition~\ref{prop:noncontracting2}, the globular subdivision $f$ induces a continuous map $\mathcal{G}^-_\alpha(X) \to \mathcal{G}^-_{f(\alpha)}(Y)$ which is a transfinite composition of homeomorphisms. Hence the proof is complete.
\epf

\begin{cor} \label{cor:enfinfin}
	Let $f:X\to Y$ be a globular subdivision. Then for any $n\geq 0$, there is the isomorphism $H_n^-(f):H_n^-(X)\iso H_n^-(Y)$.
\end{cor}

\bpf
Let $f:X\to Y$ be a globular subdivision, which means that both $X$ and $Y$ are cellular and that $|f|$ is a homeomorphism. By Theorem~\ref{thm:noncontracting3}, for any $\alpha\in X^0$, $f$ induces a homeomorphism $\mathcal{G}^-_\alpha(X)\iso \mathcal{G}^-_{f(\alpha)}(Y)$. For any $\alpha\in Y^0\backslash X^0$, the branching space $\mathcal{G}^-_\alpha(Y)$ is a singleton, $\alpha$ being inside a globular cell of $X$. We deduce that for $n\geq 1$, there are the isomorphisms 
\begin{multline*}
	H_{n+1}^-(X) = H_n(\mathcal{G}^-(X)) \iso \coprod_{\alpha\in X^0} H_n(\mathcal{G}^-_\alpha(X)) \\\iso \coprod_{\alpha\in Y^0} H_n(\mathcal{G}^-_\alpha(Y)) \iso H_n(\mathcal{G}^-(Y)) = H_{n+1}^-(Y).
\end{multline*}
For the same reason, no state $\alpha\in Y^0\backslash X^0$ is final. Thus every $\alpha \in Y^0\backslash X^0$ belongs to the image of the augmentation $\epsilon$. This implies the isomorphism $H_0^-(X)\iso H_0^-(Y)$. For $\alpha\in Y^0\backslash X^0$, $\mathcal{G}^-_\alpha(Y)$ is a singleton. Thus, the image of the map $\partial:\mathbb{Z}[\top([0,1],\mathcal{G}^-_\alpha(Y))]\to \mathbb{Z}[\mathcal{G}^-_\alpha(Y)]$ and the kernel of the map $\epsilon:\mathbb{Z}[\mathcal{G}^-_\alpha(Y)]\to \mathbb{Z}[Y^0]$ are zero for $\alpha\in Y^0\backslash X^0$. Hence we deduce that 
\[
H_1^-(X)=\coprod_{\alpha\in X^0} G^\alpha H_1^-(X) \iso \coprod_{\alpha\in X^0} G^{f(\alpha)} H_1^-(Y) \iso \coprod_{\alpha\in Y^0} G^\alpha H_1^-(Y) = H_1^-(Y).
\]
\epf


\begin{thebibliography}{10}
	
	\bibitem{CohomologySmallCategories}
	H.-J. Baues and G.~Wirsching.
	\newblock Cohomology of small categories.
	\newblock {\em J. Pure Appl. Algebra}, 38:187--211, 1985.
	\newblock \href {https://doi.org/10.1016/0022-4049(85)90008-8}
	{\path{https://doi.org/10.1016/0022-4049(85)90008-8}}.
	
	\bibitem{reversal-time-invariance}
	C.~Calk, E.~Goubault, and P.~Malbos.
	\newblock Time-reversal homotopical properties of concurrent systems.
	\newblock {\em Homology Homotopy Appl.}, 22(2):31--57, 2020.
	\newblock \href {https://doi.org/10.4310/HHA.2020.v22.n2.a2}
	{\path{https://doi.org/10.4310/HHA.2020.v22.n2.a2}}.
	
	\bibitem{Tyrone}
	T.~Cutler.
	\newblock Convergent sequence of paths in a subcomplex.
	\newblock MathOverflow.
	\newblock \url{https://mathoverflow.net/q/497175}.
	
	\bibitem{dubut_PhD}
	J.~Dubut.
	\newblock {\em Directed homotopy and homology theories for geometric models of
		true concurrency}.
	\newblock PhD thesis, Universit\'e de Saclay, Ecole Normale Sup\'erieure de
	Cachan, 2017.
	
	\bibitem{NaturalHomology}
	J.~Dubut, E.~Goubault, and J.~Goubault-Larrecq.
	\newblock Natural homology.
	\newblock In {\em Automata, languages, and programming. 42nd international
		colloquium, ICALP 2015, Kyoto, Japan, July 6--10, 2015. Proceedings. Part
		II}, pages 171--183. Berlin: Springer, 2015.
	\newblock \href {https://doi.org/10.1007/978-3-662-47666-6_14}
	{\path{https://doi.org/10.1007/978-3-662-47666-6_14}}.
	
	\bibitem{HomotopicalCategory}
	W.~G. Dwyer, P.~S. Hirschhorn, D.~M. Kan, and J.~H. Smith.
	\newblock {\em Homotopy limit functors on model categories and homotopical
		categories}, volume 113 of {\em Mathematical Surveys and Monographs}.
	\newblock American Mathematical Society, Providence, RI, 2004.
	\newblock \href {https://doi.org/10.1090/surv/113}
	{\path{https://doi.org/10.1090/surv/113}}.
	
	\bibitem{DAT_book}
	L.~Fajstrup, E.~Goubault, E.~Haucourt, S.~Mimram, and M.~Raussen.
	\newblock {\em Directed algebraic topology and concurrency. {With} a foreword
		by {Maurice} {Herlihy} and a preface by {Samuel} {Mimram}}.
	\newblock SpringerBriefs Appl. Sci. Technol. Springer, 2016.
	\newblock \href {https://doi.org/10.1007/978-3-319-15398-8}
	{\path{https://doi.org/10.1007/978-3-319-15398-8}}.
	
	\bibitem{FR}
	L.~Fajstrup and J.~Rosick{\'y}.
	\newblock A convenient category for directed homotopy.
	\newblock {\em Theory Appl. Categ.}, 21:7--20, 2008.
	
	\bibitem{model3}
	P.~Gaucher.
	\newblock A model category for the homotopy theory of concurrency.
	\newblock {\em Homology Homotopy Appl.}, 5(1):p.549--599, 2003.
	\newblock \href {https://doi.org/10.4310/hha.2003.v5.n1.a20}
	{\path{https://doi.org/10.4310/hha.2003.v5.n1.a20}}.
	
	\bibitem{exbranching}
	P.~Gaucher.
	\newblock Homological properties of non-deterministic branchings and mergings
	in higher dimensional automata.
	\newblock {\em Homology Homotopy Appl.}, 7(1):51--76, 2005.
	\newblock \href {https://doi.org/10.4310/hha.2005.v7.n1.a4}
	{\path{https://doi.org/10.4310/hha.2005.v7.n1.a4}}.
	
	\bibitem{hocont}
	P.~Gaucher.
	\newblock Inverting weak dihomotopy equivalence using homotopy continuous flow.
	\newblock {\em Theory Appl. Categ.}, 16(3):pp 59--83, 2006.
	
	\bibitem{3eme}
	P.~Gaucher.
	\newblock T-homotopy and refinement of observation ({III}): {I}nvariance of the
	branching and merging homologies.
	\newblock {\em New York J. Math.}, 12:319--348, 2006.
	
	\bibitem{mdtop}
	P.~Gaucher.
	\newblock Homotopical interpretation of globular complex by multipointed
	d-space.
	\newblock {\em Theory Appl. Categ.}, 22(22):588--621, 2009.
	
	\bibitem{leftproperflow}
	P.~Gaucher.
	\newblock Left properness of flows.
	\newblock {\em Theory Appl. Categ.}, 37(19):562--612, 2021.
	
	\bibitem{QHMmodel}
	P.~Gaucher.
	\newblock Six model categories for directed homotopy.
	\newblock {\em Categ. Gen. Algebr. Struct. Appl.}, 15(1):145--181, 2021.
	\newblock \href {https://doi.org/10.52547/cgasa.15.1.145}
	{\path{https://doi.org/10.52547/cgasa.15.1.145}}.
	
	\bibitem{Moore3}
	P.~Gaucher.
	\newblock Homotopy theory of {M}oore flows ({III}).
	\newblock {\em North-West. Eur. J. Math.}, 10:55--113, 2024.
	
	\bibitem{cubical-branching}
	P.~Gaucher.
	\newblock Branching spaces of precubical sets, 2025.
	\newblock \href {https://doi.org/10.48550/arXiv.2508.14839}
	{\path{https://doi.org/10.48550/arXiv.2508.14839}}.
	
	\bibitem{MultipointedSubdivision}
	P.~Gaucher.
	\newblock Globular subdivisions are dihomotopy equivalences.
	\newblock {\em Bull. Belg. Math. Soc. - Simon Stevin}, 32(5):607--656, 2025.
	\newblock \href {https://doi.org/10.36045/j.bbms.250522}
	{\path{https://doi.org/10.36045/j.bbms.250522}}.
	
	\bibitem{GlobularNaturalSystem}
	P.~Gaucher.
	\newblock Natural homotopy of multipointed d-spaces.
	\newblock {\em Math. J. Okayama Univ.}, 68:13--62, 2026.
	
	\bibitem{diCW}
	P.~Gaucher and E.~Goubault.
	\newblock Topological deformation of higher dimensional automata.
	\newblock {\em Homology Homotopy Appl.}, 5(2):39--82 (electronic), 2003.
	\newblock Algebraic topological methods in computer science (Stanford, CA,
	2001).
	\newblock \href {https://doi.org/10.4310/HHA.2003.v5.n2.a3}
	{\path{https://doi.org/10.4310/HHA.2003.v5.n2.a3}}.
	
	\bibitem{zbMATH07226006}
	E.~{Goubault} and S.~{Mimram}.
	\newblock {Directed homotopy in non-positively curved spaces}.
	\newblock {\em {Log. Methods Comput. Sci.}}, 16(3):55, 2020.
	\newblock Id/No 4.
	\newblock \href {https://doi.org/10.23638/LMCS-16(3:4)2020}
	{\path{https://doi.org/10.23638/LMCS-16(3:4)2020}}.
	
	\bibitem{mg}
	M.~Grandis.
	\newblock Directed homotopy theory. {I}.
	\newblock {\em Cah. Topol. G\'eom. Diff\'er. Cat\'eg.}, 44(4):281--316, 2003.
	
	\bibitem{MR99h:55031}
	M.~Hovey.
	\newblock {\em Model categories}.
	\newblock American Mathematical Society, Providence, RI, 1999.
	\newblock \href {https://doi.org/10.1090/surv/063}
	{\path{https://doi.org/10.1090/surv/063}}.
	
	\bibitem{vstrom1}
	A.~Str{\o}m.
	\newblock Note on cofibrations.
	\newblock {\em Math. Scand.}, 19:11--14, 1966.
	\newblock \href {https://doi.org/10.7146/math.scand.a-10791}
	{\path{https://doi.org/10.7146/math.scand.a-10791}}.
	
\end{thebibliography}
\end{document}